\documentclass[12pt,centertags,oneside]{amsart}
\usepackage{amsmath,amstext,amsthm,amscd,typearea,hyperref}
\usepackage{amssymb}
\usepackage{a4wide}
\usepackage[mathscr]{eucal}
\usepackage{mathrsfs}
\usepackage{typearea}
\usepackage{charter}
\usepackage{pdfsync}
\usepackage{stmaryrd}
 \usepackage[T1]{fontenc}
\usepackage{enumitem}

\usepackage{calc}

\usepackage{changepage}

\usepackage[a4paper,width=16.5cm,top=3cm,bottom=3cm]{geometry}

\numberwithin{equation}{section}

\usepackage[french,english]{babel}

\allowdisplaybreaks
\emergencystretch=\maxdimen
\usepackage{multicol}

\usepackage{xcolor}

\newtheorem{theorem}{Theorem}[section]

\newtheorem{proposition}[theorem]{Proposition}

\newtheorem{lemma}[theorem]{Lemma}

\newtheorem*{definition*}{Definition}

\newcommand{\cali}[1]{\mathscr{#1}}

\renewcommand{\Re}{\mathop{\mathrm{Re}}}

\newcommand{\dist}{\mathop{\mathrm{dist}}\nolimits}

\newcommand{\ddc}{{\rm dd^c}}
\newcommand{\dc}{{\rm d^c}}

\newcommand{\dd}{{\rm d}}
\newcommand{\Sym}{{\rm Sym}}
\newcommand{\Cyc}{{\rm Cyc}}
\newcommand{\dbar}{\overline\partial}

\newcommand{\id}{{\rm id}}

\newcommand{\ep}{\epsilon}

\newcommand{\Var}{{\rm Var}}

\newcommand{\Cc}{\cali{C}}

\newcommand{\Ec}{\cali{E}}
\newcommand{\Fc}{\cali{F}}
\newcommand{\Gc}{\cali{G}}

\newcommand{\Lc}{\cali{L}}

\newcommand{\Qc}{\cali{Q}}
\newcommand{\Rc}{\cali{R}}

\newcommand{\Tc}{\cali{T}}

\newcommand{\B}{\mathbb{B}}

\newcommand{\D}{\mathbb{D}}
\newcommand{\C}{\mathbb{C}}

\newcommand{\N}{\mathbb{N}}

\newcommand{\E}{\mathbf{E}}

\newcommand{\lp}{\langle}
\newcommand{\rp}{\rangle}

\newcommand{\norm}[1]{\lVert#1\rVert}

\newcommand{\sgn}{{\rm sgn}}

\newcommand{\oA}{\mathcal{A}}
\newcommand{\oB}{\mathcal{B}}

\newcommand{\oP}{\mathcal{P}}
\newcommand{\oN}{\mathcal{N}}

\newcommand{\oL}{\mathcal{L}}

\newcommand{\oI}{\mathcal{I}}
\newcommand{\oJ}{\mathcal{J}}

\newcommand{\fh}{\mathfrak{h}}

\newcommand{\bx}{\mathbf{x}}
\newcommand{\bE}{\mathbf{E}}

\newcommand{\tG}{\widetilde{G}}

\newcommand{\ee}{\mathsf e}
\newcommand{\fb}{\mathfrak b}



\usepackage{fancyhdr}
\title{Central limit theorem for energies of determinantal point processes on complex manifolds I}

\author{Hao Wu}
\address{School of Mathematics,  Nanjing University - Nanjing - China 210093}
\email{haowu@nju.edu.cn}

\date{}
\thanks{}

\begin{document}

\begin{abstract}
 We study the Bergman determinantal point process associated with a positive line bundle on a compact Riemann surface. Consider the U-statistic given by the Green function, called the Green energy. For a configuration of $N$ points, the expectation of Green energy is $\frac12 N\log N+O(N)$, and the variance is $(\log 2-\pi^2/24)N+O(\log^2 N)$. We prove a central limit theorem for the Green energy by combining local Bergman kernel asymptotics with an expansion for the cumulants. A central limit theorem also holds for the logarithmic energy of Bergman determinantal point process on complex manifolds.
\end{abstract}

\clearpage\maketitle
\thispagestyle{empty}


\medskip


\setcounter{tocdepth}{1}
 
  \tableofcontents\normalfont  

\section{Introduction}

The theory of determinantal point processes is a classical and well-developed topic in the past several decades. Early examples arose in random matrix theory, initiated by Ginibre \cite{Ginibre} in the study of eigenvalues of complex Gaussian matrices. The early general framework was provided by Macchi \cite{Macchi} on the fermion processes.

The determinantal point process of a given space $X$ is defined by a probability measure on the space $X^N$,  the configuration space of $N$ particles on $X$, whose density can be written as
\begin{align*}
\rho^{(N)} (x_1,x_2,\dots,x_N)={1\over N!} \det \big(K_N(x_j,x_k) \big)_{1\leq j,k\leq N},
\end{align*}
where $K_N$ is the integral kernel of an orthogonal projection operator onto a vector space of dimension $N$. The reader may refer to the surveys \cite{hou-kri-man,Soshnikov} for more background.

A classical example is the complex Ginibre ensemble. Let $z_1,\ldots,z_N$ be the eigenvalues, listed in uniformly random order, of an $N\times N$ random matrix whose entries are independent centered complex Gaussian variables with variance $1/N$. With respect to the  Lebesgue measure  on $\C^N$, their joint density is
$$\rho^{(N)}(z_1,\ldots,z_N)
=\frac{1}{N!}\det\!\big(K_N(z_j,z_k)\big)_{1\leq j,k\leq N}
=\frac{N^{N(N+1)/2}}{\pi^N\prod_{j=1}^{N}j!}\,
e^{-N\sum_{k=1}^{N}|z_k|^2}
\prod_{1\leq i<j\leq N}|z_j-z_k|^2,$$
where
$$K_N(z,w)=\frac{N}{\pi}\,e^{-\frac N2(|z|^2+|w|^2)}
\sum_{k=0}^{N-1}\frac{(Nz\overline w)^k}{k!}.$$
 Thus the eigenvalues form a determinantal point process with kernel \(K_N\).

\medskip

In complex geometry, Berman \cite{ber-cmp2} studied the \textit{Bergman determinantal point processes} on a compact complex manifold, whose correlation kernels are the Bergman kernels associated with powers of a positive holomorphic line bundle. He  established  Gaussian fluctuations for linear statistics, base on local estimates of the Bergman kernels.

In this paper, we consider the Green energy of the Bergman determinantal point process on a compact  Riemann surface,  which is a second-order U-statistic obtained by summing the Green functions over all ordered pairs of distinct particles. We compute the leading coefficient of its variance and prove that the centered energy, divided by its standard deviation, converges in distribution to a standard Gaussian. We will also give some remarks for the situation in higher dimensions. 

\medskip
 
Let $X$ be a compact Riemann surface and let $\oL$ be a positive holomorphic line bundle on $X$ with $\deg(\oL) \geq 1$.  Fix a Hermitian metric $\fh$ on $\oL$ such that the curvature form $c_1(\oL,\fh)$ is strictly positive. The normalized $(1,1)$-form
$$ \omega :=c_1(\oL,\fh)/\deg(\oL)$$
is a smooth probability measure on $X$.

 For each $n\in\N$, consider the space of global holomorphic sections $H^0 (X,\oL^n)$ of $\oL^n$.  Its dimension  is
 \[
 N_n:=\dim_{\mathbb{C}} H^0(X, \oL^n) = n \cdot \deg(\oL) - g + 1
 \] 
 by Riemann-Roch theorem, where  $g$ is the genus of $X$.
 For each element $s\in  H^0(X, \oL^n)$, the norm at each point $x\in X$, induced by $\fh$, is denoted by $|s(x)|_\fh$. The line bundle $\oL^n$ is also positive. Its degree is equal to $n\deg (\oL)$. The metric $\fh$ on $\oL$ induces a metric $\fh_n$ on $\oL^n$ whose curvature form is equal to $n\deg(\oL)\omega$. It is defined by $\norm{s^{n}(x)}_{\fh_n}:=\norm{s(x)}_\fh^n$ for any section $s$ of $\oL$. We use the notation $(\cdot,\cdot)_n$ to denote the hermitian inner product at each point $x\in X$, which is given by $\fh_n$.  The metric $\fh$, together with $\omega$, induce an inner product $\lp\cdot,\cdot\rp_n$ on the complex vector space $H^0(X,\oL^n)$ of holomorphic sections of $\oL^n$:
 \begin{equation*}
 \lp s_1,s_2\rp_n:=\int_{x\in X} (s_1(x),s_2(x))_n \,\omega(x)     \quad  \text{for} \quad  s_1,s_2\in H^0(X,\oL^n).
 \end{equation*}
Let $S_n:=(s_1,s_2,\dots,s_{N_n})$  be an orthonormal  basis  of $H^0 (X,\oL^n)$. 
 For any points $x_1,x_2,\dots,x_{N_n}\in X$, define
$$  \det S_n (x_1,x_2,\dots,x_{N_n}):= \det (s_j (x_\ell))_{1\leq j,\ell\leq N_n},$$
which is a section of $(\oL^n)^{\boxtimes N_n}$ over $X^{N_n}$. The metric $\fh$ induces in a canonical way a metric  $(\fh^n)^{\boxtimes N_n}$ on $(\oL^n)^{\boxtimes N_n}$. Denote this norm by $\big|\det S_n(x_1,\dots,x_{N_n})\big|_\fh$.

 Consider the point process on $X$, which is induced by the following probability measure on $X^{N_n}$:
\begin{align*}\label{correla}
\rho^{N_n}(x_1,x_2,\dots,x_{N_n}):= \mathcal Z_{N_n}  \big|\det S_n(x_1,\dots,x_{N_n})\big|_\fh^2 \,\omega(x_1)\omega(x_2)\cdots\omega(x_{N_n}) 
\end{align*}    
for some normalization constant $\mathcal Z_{N_n}>0$. One will see later in \eqref{kenel-density} that this process is determinantal.

For   $\bx:=(x_1,x_2,\dots,x_{N_n})\in X^{N_n}$, we define the \textit{Green energy} of $\bx$ by
$$  \Ec_{N_n}(\mathbf x):= \sum_{1\leq j,k\leq N_n, j\neq k}G (x_j,x_k),   $$
where $G$ is the \textit{Green function} of $(X,\omega)$, see  \eqref{Green}. Each unordered pair occurs twice in this convention.
It can also be  viewed as a random variable for the  determinantal  point process defined above. The Green energy plays an important role in the theory of Bergman determinantal point process, as it relates to the rate function in the large deviation principle \cite{ber-cmp}. It  also appears as the leading term in the bosonization formula \cite{Alv-Luis,kle-cmp,wu-grow}.

To simplify the notation, we will abbreviate $N_n$ by $N$ throughout this article.   Our main result is the following central limit theorem.

\begin{theorem}\label{thm-main}
	As $N\to \infty$,  
	$$ { \Ec_N- \bE(\Ec_N) \over \sqrt{ \Var (\Ec_N)    } } \longrightarrow  \mathcal N(0,1)  \quad\text{in law},  $$
	with
	$$\bE(\Ec_N)=\frac12N\log N+O(N)$$
    and 
   	$$\Var(\Ec_N)=\Big(\log 2-{\pi^2\over 24}\Big)N+O(\log^2 N).$$
\end{theorem}

In the proof of the CLT in \cite{ber-cmp2}, the local estimate of Bergman kernel in a window of size $1/\sqrt N$  play an essential role. The application of Bergman kernel can be found in \cite{auv-ma-mann,ber-bou-invent,ber-acta,ber-ajm,cat-ber,csm,co-lu-adv,dai-liu-ma,dinh-jfa1,din-ma-ngu-ens,fin1,fin-mann,ma-mar-mann,mar-sav-mann,mar-vu-crelle,tian-jdg,vu-amj,zel-imrn-sze} and the standard book by Ma-Marinescu \cite{ma-book}.
However,  the approach for linear statistics cannot be applied to U-statistics. The expectation and variance for linear statistics can be easily obtained from the local estimates of Bergman in a window of size $1/\sqrt N$. While, for the Green energy, there is  a logarithmic singularity when two particles meet each other. The window of size $1/\sqrt N$ is not enough to bound the singular integral. Instead, we will use an estimate of Bergman kernel in a window of size $\sqrt{\log N/N}$ established by Shiffman-Zelditch \cite{shi-zel-gafa} and Lu-Shiffman \cite{lu-shi-jga}. Besides, Berman proved the central limit theorem by building a clever and quick method to show the cumulant generating function converges to $t^2/2$, which avoids the long computation of the higher moments and cumulants of the linear statistic. In the case of U-statistics, the method fails. To control the higher cumulants of the Green energy, we have to write down the expansion of higher cumulants of the Green energy term by term, and prove an upper bound for it using some combinatorial techniques from Krishnapur-Yogeshwaran \cite{man-mom}.

Unfortunately, our proof does not apply to the case when the metric $\fh$ on $\oL$ is singular, because the smoothness condition of the Green function outside the diagonal is necessary.

\medskip

This work was inspired by Michelen-Yakir \cite{Mich-Yak}, who studied
 the logarithmic energy   of the zeros
of SU($2$) Gaussian polynomials on the Riemann sphere. They proved a CLT and obtained asymptotic formulas for the expectation and variance.    The point process induced by the zeros of random polynomials or random holomorphic sections on complex manifolds has been studied extensively in last twenty years, see e.g.\ \cite{bay-com-mar-TAMS,buc-nis-ron-sod-ptrf,coman-marin-Nguyen,DGV-2024,dinh-jfa1,din-mar-sch-jsp,din-sib-cmh,dlm-jussieu,gho-nis-cpam,Nishry-Wennman,shi-zel-cmp,shi-zel-gafa,shi-zel-zre-ind,sod-tsi-I,zei-zel-imrn,zel-imrn} and the surveys \cite{survey-random-hol,gaussian-book,geo-duc-sur}.

Central limit theorems for U-statistics of determinantal point
processes have been established in several settings. Aguirre-Soshnikov-Sumpter \cite{agu-sos-sum} proved  the CLT for smooth pair statistics of
the circular unitary ensemble. Later, Soshnikov-Wu \cite{sos-wu} obtained Gaussian
fluctuations for  pair-counting statistics of circular unitary ensemble.
Feng-G\"otze-Yao \cite{feng-goe-yao} proved a CLT
for bounded symmetric pair kernels for spherical projection processes. 

There are many works of CLT for linear statistics of determinantal point
processes. We only list a few of them here 
\cite{ame-hed,bai-wan-1,bai-wan-2,ber-toep,cos-leb,forpj,joh,leb-ser,li-xu,rid-vir,sos-clt}. See also the survey \cite{hou-kri-man} and the references therein.

\medskip

When $(X,\omega)$ is a K\"ahler manifold of dimension $d\geq 1$, the Green function \eqref{Green} is not defined. Instead, for $\bx:=(x_1,x_2,\dots,x_{N})\in X^{N}$, we may consider the \textit{logarithmic energy}  
$$  \Ec^{\log}_{N}(\mathbf x):= \sum_{1\leq j,k\leq N, j\neq k}\log\dist (x_j,x_k),   $$
where $\dist$ is the distance function induced by the K\"ahler form $\omega$.    A similar central limit theorem holds for this energy. To keep the present paper in a reasonable length, we postpone it to the forthcoming work \cite{wu-clt-2} and give a short remark in Section \ref{sec-remk}.

\medskip

\noindent\textbf{Some notation.}
The symbols $\lesssim$ and $\gtrsim$ stand for inequalities up to a positive multiplicative constant. If both hold, then we write $\simeq$. We write 
$A_N \overset{C_N}{\asymp}  B_N$ if $A_N$ and $B_N$ differ by a term $O(C_N)$ as $N\to\infty$.   The dependence of this constant on certain parameters, or lack thereof, will be clear from the context. We denote by $\D$ the unit disc in $\C$ and $\D(a,r)$ the open disc with center $a$ and radius $r$.  For other Riemann surfaces or complex manifolds, $\B(a,r)$ denotes the ball of center $a$ and radius $r$.  The expectation of a random variable or a random measure is denoted by $\E(\cdot)$.  Recall that $\dc:=\frac{i}{2\pi}(\dbar-\partial)$ and $\ddc=\frac{i}{\pi}\partial\dbar$.

\medskip
\noindent\textbf{Acknowledgements.}  The author would like to thank Professor Duc-Viet Vu for his valuable remark on  Toeplitz operator, which accelerated the calculation of variance.

\medskip
\noindent\textbf{AI statement.} This work is independent of AI.

\medskip
\noindent\textbf{Funding.}
This work was funded by National Key R\&D Program of China 2025YFA1018300 and Basic Research Program of Jiangsu BK20251226.

\section{Auxiliary integrals for  the Bergman kernel}

In this section, we do all the detailed computations   related to the Bergman kernel in a local normal coordinate needed in Section \ref{sec-pre}.

Define a function $E_1$ on $(0,\infty)$ and a constant $\gamma_E$ by 
$$  E_1(a)=\int_{a}^\infty t^{-1} e^{-t} \,\dd t, \quad \gamma_E:=-\int_0^\infty e^{-t}\log t\,\dd t.$$
Note that $E_1(0)$ is not defined and as $a\to 0^+$,
\begin{align*}
E_1(a)=-\log a-\gamma_E +O(a).
\end{align*}

\begin{lemma}\label{lem-E1}
	For $a>0$, we have
	$$ 2\pi N^2 \int_a ^\infty re^{-\pi N r^2} \log r \,\dd r=Ne^{-\pi Na^2}\log a+{N\over 2}E_1(\pi Na^2)  $$
	and
	$$ 2\pi N^2 \int_0 ^\infty re^{-\pi N r^2} \log r \,\dd r=-{N\log N\over 2}-{\log \pi +\gamma_E\over 2}N.$$
\end{lemma}

\begin{proof}
	Substituting $t=\pi N r^2$, the left hand side becomes
	\begin{align*}
	N\int_{\pi Na^2}^\infty e^{-t} &\log\sqrt{t\over \pi N} \,\dd t 
	={N\over 2}\int_{\pi Na^2}^\infty  e^{-t} \log t  \,\dd t-{N\over 2}\log{\pi N}\int_{\pi Na^2}^\infty  e^{-t}   \,\dd t \\
	&={N\over 2} \big[ -e^{-t}\log t\big]_{t=\pi N a^2}^\infty +{N\over 2} \int_{\pi Na^2}^\infty  e^{-t} t^{-1} \,\dd t -{N\over 2}\log(\pi N) e^{-\pi Na^2} \\
	&={N\over 2} \log(\pi N a^2) e^{-\pi Na^2}+{N\over 2}E_1(\pi Na^2)-{N\over 2}\log(\pi N) e^{-\pi Na^2} \\
	&=Ne^{-\pi Na^2}\log a+{N\over 2}E_1(\pi Na^2).
	\end{align*}
	
	When $a=0$,
	\begin{align*}
	2N^2 \pi\int_0 ^{\infty}  r  e^{-\pi N  r^2}   \log r   \,\dd r&={N\over 2}\int_{0}^\infty  e^{-t} \log t  \,\dd t-{N\over 2}\log{\pi N}\int_{0}^\infty  e^{-t}   \,\dd t\\
	&=-{N\over 2}\gamma_E-{N\over 2}\log(\pi N)
	=-{N\log N\over 2}-{\log \pi +\gamma_E\over 2}N.
	\end{align*}
	This ends the proof.
\end{proof}

For $x\in X$, we say  $z$ is a local normal  coordinate centering at $x$ with respect to  $\omega$ if   
\begin{equation}\label{defn-normal}
z(x)=0 \quad\text{and}\quad \omega \big(z(w)\big) = {i\over 2}\big(1+O(|z(w)|^2)\big) \,\dd z(w) \wedge  \dd \overline{z(w)} 
\end{equation}
for $w$ near $x$.
Such a coordinate always exists on K\"ahler manifolds.

\begin{lemma}\label{lem-intE1}
	Let $z$ be a local normal coordinate centering at $x$ containing $\B(x,r_0)$ for a fixed $r_0>0$.
	We have
	$$\int_{|z|<r_0}   E_1(\pi N|z|^2)\,\omega(z) =N^{-1}+O(N^{-2}),    $$
	$$\int_{|z|<r_0} \big| \log |z| \big| E_1(\pi N|z|^2)\,\omega(z) =O(N^{-1}\log N),    $$
	and
	$$\int_{|z|<r_0} \Big(  E_1(\pi N|z|^2) \Big)^2\,\omega(z)={ 2\log 2 \over N}+O(N^{-2}).$$
\end{lemma}

\begin{proof}
	By definition,
	\begin{align*}
	\int_{|z|<r_0}   E_1(\pi N|z|^2)\,\omega(z) &=\int_{|z|<r_0}  \int_{\pi N|z|^2}^\infty  t^{-1}e^{-t}\,\dd t \,{i\over 2}\big(1+O(|z|^2)\big)\dd z\wedge \dd\overline z \\
	&=\int_0^{2\pi}\int_0^{r_0}  \int_{\pi Nr^2}^\infty t^{-1}e^{-t}\,\dd t \, \big(1+O(r^2)\big) r\,\dd r\dd\theta \\
	&=2\pi \int_0^{\pi Nr_0^2} {1\over 2\pi N} \int_{u}^\infty t^{-1}e^{-t}  \big(1+O(N^{-1}u)\big) \,\dd  t \dd u \\
	&\overset{N^{-2}}{\asymp}  {1\over  N} \int_0^{\infty}  \int_{u}^\infty t^{-1}e^{-t}  \big(1+O(N^{-1}u)\big) \,\dd  t \dd u \\
	&={1\over N}\int_{0}^{\infty} \int_0^t  t^{-1}e^{-t}  \big(1+O(N^{-1}u)\big)  \,\dd u\dd t  \\
	&={1\over N}\int_{0}^{\infty} \big(e^{-t}+O(N^{-1}te^{-t}/2)\big)\,\dd t \\&=N^{-1}+O(N^{-2}).
	\end{align*}
	This proves the first identity. The second one is similar by adding a $\log |z|$ factor in last computation, which produces an extra $\log N$ factor for the integral.  For the last one,
	\begin{align*}
	&\int_{|w|<r_0} \Big(  E_1(\pi N|w|^2) \Big)^2\,\omega(w)= \int_{|w|<r_0} \int_{\pi N|w|^2}^\infty\int_{\pi N|w|^2}^\infty t^{-1} e^{-t}  s^{-1}e^{-s}\,\dd t \dd s \,\omega(w)\\
	&=\int_0^{2\pi}\int_0^{r_0}   \int_{\pi Nr^2}^\infty\int_{\pi Nr^2}^\infty t^{-1} e^{-t}  s^{-1}e^{-s}\,\dd t \dd s \big(1+O(r^2)\big) r \, \dd r \dd\theta \\
	&=2\pi \int_0^{\pi Nr_0^2} {1\over 2\pi N}\int_{u}^\infty  \int_{u}^\infty t^{-1} e^{-t}  s^{-1}e^{-s}\big(1+O(N^{-1}u)\big)\,\dd t \dd s   \dd u \\
	&\overset{N^{-2}}{\asymp}  {1\over N} \int_0 ^\infty \int_{u}^\infty  \int_{u}^\infty t^{-1} e^{-t}  s^{-1}e^{-s}\big(1+O(N^{-1}u)\big)\,\dd t \dd s   \dd u\\
	&= {1\over N} \int_0 ^\infty \int_{u}^\infty  \int_{u}^\infty t^{-1} e^{-t}  s^{-1}e^{-s} \,\dd t \dd s   \dd u+O(N^{-2}).
	\end{align*}
	This integral is symmetric in $t$ and $s$. Thus, it is double of the integral over the sub-domain $\{t\geq s\}$, which is
	\begin{align*}
	{2\over N} \int_0 ^\infty \int_{u}^\infty  \int_{s}^\infty t^{-1} e^{-t}  s^{-1}e^{-s} \,\dd t \dd s   \dd u &={2\over N}  \int_0 ^\infty \int_{0}^t  \int_{0}^s t^{-1} e^{-t}  s^{-1}e^{-s} \,\dd u \dd s   \dd t \\
	={2\over N}  \int_0 ^\infty \int_{0}^t t^{-1} e^{-t}e^{-s} \,  \dd s   \dd t 
	&={2\over N}  \int_0 ^\infty   t^{-1} e^{-t} (1-e^{-t})      \, \dd t \\
	={2\over N}  \int_0 ^\infty t^{-1}(e^{-t}-e^{-2t})\,\dd t&={2\log 2\over N} 
	\end{align*}
	showing the last identity.
\end{proof}

We need the following formula for radial functions.

\begin{lemma}\label{lem-radial}
	Let $\varphi(x,y)$ be a  radial function on $\C\times \C$, i.e., $\varphi(x,y)=\varphi(|x|,|y|)$ for $x,y\in \C$. Suppose $\varphi(x,y) e^{|xy|}$ is integrable. Then 
	$$\int \varphi(x,y) e^{x\overline y} \,{i\over 2}\dd x\wedge \dd\overline x\,{i\over 2}\dd y\wedge \dd\overline y=4\pi^2\int_0^\infty\int_0^\infty\varphi(s,t )st \, \dd s \dd t.$$
	If $\phi(x,y)$ is radial on $x$, i.e., $\phi(x,y)=\phi(|x|,y)$ for $x,y\in\C$, and $\phi(x,y_0)e^{|x  y_0|}$ is integrable  for some fixed $y_0$, then
	$$\int \phi(x,y_0)e^{x  y_0} \,{i\over 2}\dd x\wedge \dd\overline x
	=2\pi \int_0^\infty \phi(s,y_0) s\,\dd s.$$
\end{lemma}

\begin{proof}
	Using polar coordinates, we have
	\begin{align*}
	\int \varphi(x,y) e^{x\overline y} \,{i\over 2}\dd x\wedge \dd\overline x\,{i\over 2}\dd y\wedge \dd\overline y =\int_0^{2\pi}\int_0^{2\pi}\int_0^\infty\int_0^\infty\varphi(|r_x|,|r_y| )e^{r_xr_ye^{i\theta_x-i\theta_y}}  r_xr_y \,\dd r_x \dd r_y\dd\theta_x \dd \theta_y\\
	=\int_0^\infty\int_0^\infty\varphi(|r_x|,|r_y| )r_xr_y \Big( \int_0^{2\pi}\int_0^{2\pi} e^{r_xr_ye^{i\theta_x-i\theta_y}}\,\dd\theta_x \dd \theta_y \Big)\dd r_x \dd r_y.
	\end{align*}
	Write the function 
	$$e^{r_xr_ye^{i\theta_x-i\theta_y}} =\sum_{k=0}^\infty  r_x^kr_y^{k} { e^{ik\theta_x-ik\theta_y} \over k! }.  $$
	Notice that for all $k\neq 0$,
	$$ \int_0^{2\pi}\int_0^{2\pi} r_x^kr_y^k { e^{ik\theta_x-ik\theta_y} \over k! }\,\dd\theta_x \dd \theta_y  =0$$ 
	Thus, the integral of $\varphi(x,y)$ becomes
	$$  4\pi^2\int_0^\infty\int_0^\infty\varphi(|r_x|,|r_y| )r_xr_y \, \dd r_x \dd r_y.  $$
	This gives the first equality. The second one is similar.
\end{proof}

\begin{lemma}\label{lem-ywlog}
	For $a,b,c,d,s,t\in \N$, and $N\geq 2$, we have
	\begin{align*}
	\int |y|^{a} |w|^b \big|\log|y|\big|^c
	\big|\log|w|\big|^d 
	N^3 \big|e^{-{\pi} N( |y|^2+|w|^2-y\overline  w)}\big|
	\,{i\over 2} \,\dd y \wedge  \dd \overline y\,{i\over 2}  \,\dd w \wedge  \dd \overline w\\
	\lesssim _{a,b,c,d}N^{1-a/2-b/2}(\log N)^{c+d}
	\end{align*}
	and
	\begin{align*}
	\int |y|^{a} |w|^b |z|^s &  \big|  \log|y|\big|^c
	\big|\log|w|\big|^d   \big| \log|z|\big|^t
	N^4 \big|e^{-{\pi} N( |y|^2+|w|^2+|z|^2-y\overline  w-w\overline z)}\big|\\
	&\,{i\over 2} \,\dd y \wedge  \dd \overline y\,{i\over 2}  \,\dd w \wedge  \dd \overline w\,{i\over 2}  \,\dd z \wedge  \dd \overline z 
	\lesssim _{a,b,c,d,s,t}N^{1-a/2-b/2-s/2}(\log N)^{c+d+t}.
	\end{align*}
\end{lemma}

\begin{proof}
	We only prove the first assertion. The second one can be obtained in a similar way. Using that $\Re\big(|y|^2+|w|^2-y\overline  w\big) \geq 1/2\big(|y|^2+|w|^2\big) $, the integral is bounded from above by
	\begin{align*}
	&	\int |y|^{a} |w|^b \big|\log|y|\big|^c
	\big|\log|w|\big|^d 
	N^3 e^{-{\pi N \over 2}( |y|^2+|w|^2)}
	\,{i\over 2} \,\dd y \wedge  \dd \overline y\,{i\over 2}  \,\dd w \wedge  \dd \overline w \\
	&=	\int \Big|{y\over \sqrt N} \Big|^{a} \Big|{w\over \sqrt N}\Big|^b \Big|\log {|y|\over \sqrt N} \Big| ^c
	\Big|\log {|w|\over \sqrt N} \Big| ^d 
	N e^{-{\pi  \over 2}( |y|^2+|w|^2)}
	\,{i\over 2} \,\dd y \wedge  \dd \overline y\,{i\over 2}  \,\dd w \wedge  \dd \overline w \\
	&\lesssim_{c,d}N^{1-a/2-b/2} 	\int |y|^{a} |w|^b \Big(\big|\log |y|\big|^c+\log^c\sqrt N \Big)
	\Big(\big|\log |w|\big|^d+\log^d \sqrt N \Big)
	\\
	&\qquad\qquad\qquad\qquad\qquad\qquad\qquad e^{-{\pi  \over 2}( |y|^2+|w|^2)} {i\over 2} \,\dd y \wedge  \dd \overline y\,{i\over 2}  \,\dd w \wedge  \dd \overline w. 
	\end{align*}
	We need to show that last integral is $\lesssim(\log N)^{c+d}$.  Indeed, for any $0\leq c'\leq c,0\leq d'\leq d$, the function $|y|^a|w|^b ( \log|y|)^{c'}(\log|w|)^{d'}e^{-{\pi  \over 2}( |y|^2+|w|^2)}$ is integrable on $\C^2$.  After expanding the integrand, we see the exponent of $\log  N$ is $c+d$. 
\end{proof}

\begin{lemma}\label{lem-logw-logy}
	The integral
	$$
	\int  \log|y|
	\big(\log|w|-\log|y| \big)
	N^3 e^{-{\pi} N( |y|^2+|w|^2-y\overline  w)}
	\,{i\over 2} \,\dd y \wedge  \dd \overline y\,{i\over 2}  \,\dd w \wedge  \dd \overline w=-{\pi^2\over 24}N.$$
\end{lemma}

\begin{proof}
	By Lemma \ref{lem-radial}, the integral becomes
	\begin{align*}
	4\pi^2\int_{0}^\infty  \int_{0}^\infty \log|s|(\log|t|-\log|s|) N^3e^{-\pi N(s^2+t^2)} st\,\dd s\dd t.
	\end{align*}
	Scaling \(u=\sqrt{\pi N}\,s\), \(v=\sqrt{\pi N}\,t\), we get
	\[
	4N\int_0^\infty\int_0^\infty
	uv(\log u-\log\sqrt{\pi N})(\log v-\log u)e^{-u^2-v^2}\,\dd u\dd v.
	\]
	By symmetry, it is equal to
	\[
	4N\left[
	\Big(\int_0^\infty u e^{-u^2}\log u\,du\Big)^2
	-\Big(\int_0^\infty u e^{-u^2}\log^2u\,du\Big)
	\Big(\int_0^\infty u e^{-u^2}\,du\Big)
	\right].
	\]
	Using integration by parts, one can compute
	\[
	\int_0^\infty u e^{-u^2}\log u\,du=-\frac{\gamma_E}{4},\quad 
	\int_0^\infty u e^{-u^2}\log^2u\,du
	={1\over 8}\Big(\gamma_E^2+{\pi^2\over 6}\Big),\quad \int_0^\infty u e^{-u^2}\,du={1\over 2}.
	\]
	Therefore, the result follows.
\end{proof}

\begin{lemma}\label{lem-zbarz}
	Let $a>0$ and $b,c\in\C$, then
	$$ \int e^{-a|z|^2+bz+c\overline z}\,{i\over 2}\dd z\wedge \dd \overline z={\pi\over a}e^{bc/a}$$
	and
	$$  \int \log|z| e^{-a|z|^2+bz+c\overline z}\,{i\over 2}\dd z\wedge \dd \overline z= {\pi\over 2a}\sum_{k=1}^\infty{H_k\over k!} \Big({bc\over a } \Big)^k-{\pi\over 2a}(\gamma_E+\log a)e^{bc/a},$$
	where $H_k:=\sum_{j=1}^k  1/j$ is the $k$-th harmonic number.
\end{lemma}

\begin{proof}
	Using polar coordinates, the first integral becomes
	\begin{align*}
	\int_0^{2\pi}\int_0^\infty e^{-a t^2 +bte^{i\theta}+cte^{-i\theta}} t\,\dd t \dd\theta=\int_0^\infty te^{-a t^2}\int_0^{2\pi}e^{bte^{i\theta}+cte^{-i\theta}} \, \dd\theta\dd t
	\end{align*}
	Write 
	$$e^{bte^{i\theta}+cte^{-i\theta}}= \sum_{m=0}^\infty {(bte^{i\theta}+cte^{-i\theta})^m \over m!}=\sum_{m=0}^\infty{t^m\over m!} \sum_{k=0}^m \binom{m}{k} b^k c^{m-k}e^{i(2k-m)\theta}.  $$
	Using that $$\int e^{i(2k-m)\theta}\,\dd\theta=0 \quad \text{for all}\quad 2k-m\neq 0,$$ we see that the first integral is equal to
	\begin{align*}
	&\int_0^\infty te^{-a t^2} 2\pi\sum_{k=0}^\infty {t^{2k}\over (2k)!}\binom{2k}{k}b^kc^k \,\dd t=\int_0^\infty te^{-a t^2} 2\pi\sum_{k=0}^\infty {t^{2k}\over (k)!(k)!}  b^kc^k \,\dd t\\
	&=2\pi \sum_{k=0}^\infty {b^kc^k\over (k)!(k)!}\int_0^\infty t^{2k+1} e^{-at^2}\,\dd t=2\pi \sum_{k=0}^\infty {b^kc^k\over (k)!(k)!}{k!\over 2a^{k+1}}={\pi\over a}\sum_{k=0}^\infty {1\over k!}\Big({b c\over a}\Big)^k={\pi\over a}e^{bc/a}.
	\end{align*}
	
	The second integral has one more $\log|z|$ factor. Repeating the  calculation above with this extra factor and setting  $H_0=0$, we obtain
	\begin{align*}
	2\pi \sum_{k=0}^\infty {b^kc^k\over (k)!(k)!}\int_0^\infty t^{2k+1} e^{-at^2}\log t \,\dd t
	&=  2\pi \sum_{k=0}^\infty {b^kc^k\over (k)!(k)!}{k!\over 4a^{k+1}}(H_k-\gamma_E-\log a)\\
	={\pi\over 2a}\sum_{k=0}^\infty{1\over k!} \Big({bc\over a } \Big)^k(H_k-\gamma_E-\log a)
	&={\pi\over 2a}\sum_{k=1}^\infty{H_k\over k!} \Big({bc\over a } \Big)^k-{\pi\over 2a}(\gamma_E+\log a)e^{bc/a}.
	\end{align*}
	The proof is complete.
\end{proof}

\begin{lemma}\label{lem-logw-logwz}
	The integral
	\begin{align*}
	\int\log|y| \big(\log|w-z|-\log|w|\big)   N^4 e^{-\pi N (|y|^2+|w|^2+|z|^2-w\overline y-y\overline z)}  \\ 
	\,{i\over 2}  \,\dd y \wedge  \dd \overline y \,{i\over 2}  \,\dd w \wedge  \dd \overline w\,{i\over 2}  \,\dd z \wedge  \dd \overline z=-{\pi^2\over 48} N.
	\end{align*}
\end{lemma}

\begin{proof}
	We integrate $y$ first using Lemma \ref{lem-zbarz}, getting
	\begin{align*}
	\int\log|y| e^{-\pi N (|y|^2-w\overline y-y\overline z)}   
	\,{i\over 2}  \,\dd y \wedge  \dd \overline y ={\pi\over 2\pi N}\sum_{k=1}^\infty{H_k\over k!} ({\pi Nw\overline z }) ^k-{\pi\over 2\pi N}\big(\gamma_E+\log (\pi N)\big)e^{\pi Nw\overline z}.
	\end{align*}
	
	By symmetry, 
	\begin{align*}
	\int (w\overline z)^k \log|w| N^4 e^{-\pi N ( |w|^2+|z|^2 )}\,{i\over 2}  \,\dd w \wedge  \dd \overline w\,{i\over 2}  \,\dd z \wedge  \dd \overline z=0
	\end{align*}
	for all $k>0$. Moreover,
	by \eqref{logxy-logx}, 
	\begin{align*}
	&\int e^{\pi Nw\overline z}(\log|w-z|-\log|w|) N^4 e^{-\pi N ( |w|^2+|z|^2 )}\,{i\over 2}  \,\dd w \wedge  \dd \overline w\,{i\over 2}  \,\dd z \wedge  \dd \overline z\\
	&=\int (\log|w-z|-\log|w|) N^4 e^{-\pi N ( |w|^2+|z|^2 -w\overline z)}\,{i\over 2}  \,\dd w \wedge  \dd \overline w\,{i\over 2}  \,\dd z \wedge  \dd \overline z=0.
	\end{align*}
	Hence the only thing left is 
	\begin{align*}
	\int {\pi\over 2\pi N}\sum_{k=1}^\infty{H_k\over k!} ({\pi Nw\overline z })^k \log|w-z| N^4 e^{-\pi N ( |w|^2+|z|^2 )}\,{i\over 2}  \,\dd w \wedge  \dd \overline w\,{i\over 2}  \,\dd z \wedge  \dd \overline z \\
	={1\over 2} \sum_{k=1}^\infty {H_k\over k!} \int ({\pi Nw\overline z })^k \log|w-z| N^3 e^{-\pi N ( |w|^2+|z|^2 )}\,{i\over 2}  \,\dd w \wedge  \dd \overline w\,{i\over 2}  \,\dd z \wedge  \dd \overline z.  
	\end{align*}
	By Lemma \ref{lem-wzlogwz}, last sum is
	$$-{1\over 2} \sum_{k=1}^\infty {H_k\over k!}\cdot {N (k-1)!\over 2^{k+1}}= -{N\over 4} \sum_{k=1}^\infty {H_k\over k 2^k}=-{N\over 4} \cdot {\pi^2\over 12}=-{\pi^2\over 48}N, $$ 
	where the sum $\sum_{k=1}^\infty {H_k\over k 2^k}$ is calculated by using the expansion $$ -{\log(1-x)\over (1-x)}= \sum_{k=1}^\infty H_k x^k \quad\text{for}\quad |x|<1. $$
	This ends the proof of the lemma.
\end{proof}

\begin{lemma}\label{lem-wzlogwz}
	For $k\geq 1$,
	$$\int ({\pi Nw\overline z })^k \log|w-z| N^3 e^{-\pi N ( |w|^2+|z|^2 )}\,{i\over 2}  \,\dd w \wedge  \dd \overline w\,{i\over 2}  \,\dd z \wedge  \dd \overline z=-{N (k-1)!\over 2^{k+1}}.$$
\end{lemma}

\begin{proof}
	Scaling $u=\sqrt{\pi N} w,v=\sqrt{\pi N} z$, the integral becomes
	\begin{align*}
	{N\over\pi^2}\int (  u\overline v  )^k \log{|u-v|\over \sqrt{\pi N}}  e^{-  ( |u|^2+|v|^2 )}\,{i\over 2}  \,\dd u \wedge  \dd \overline u\,{i\over 2}  \,\dd v \wedge  \dd \overline v.
	\end{align*}
	By symmetry, since $k>0$,
	$${N\over\pi^2}	\int (  u\overline v  )^k \log \sqrt{\pi N}   e^{-  ( |u|^2+|v|^2 )}\,{i\over 2}  \,\dd u \wedge  \dd \overline u\,{i\over 2}  \,\dd v \wedge  \dd \overline v=0.$$
	Thus, we only need to compute 
	\begin{align*}
	{N\over\pi^2}\int (  u\overline v  )^k \log|u-v|  e^{-  ( |u|^2+|v|^2 )}\,{i\over 2}  \,\dd u \wedge  \dd \overline u\,{i\over 2}  \,\dd v \wedge  \dd \overline v
	\end{align*}

	Passing to polar coordinates, it becomes
	\begin{align}
	&{N\over\pi^2}\int_0^\infty\int_0^\infty\int_0^{2\pi}\int_0^{2\pi}
	s^{k+1}t^{k+1}e^{ik(\theta-\phi)}
	\log|se^{i\theta}-te^{i\phi}|\,
	e^{-s^2-t^2}
	\,\dd\theta\dd\phi\dd s\dd t  \nonumber \\
	&={2N\over\pi^2}\int_0^\infty\int_0^s\int_0^{2\pi}\int_0^{2\pi}
	s^{k+1}t^{k+1}e^{ik(\theta-\phi)}
	\log|se^{i\theta}-te^{i\phi}|\,
	e^{-s^2-t^2}
	\,\dd\theta\dd\phi\dd t\dd s,
	\label{eq:polar}
	\end{align}
	where we use the fact that the integral is symmetric on $s,t$.
	Using that $|se^{i\theta}-te^{i\phi}|^2=s^2+t^2-2st\cos(\theta-\phi)$.
	when $s>t$, we have
	\begin{align*}
	\log|se^{i\theta}-te^{i\phi}|
	&=
	\frac12\log\big(s^2+t^2-2st\cos(\theta-\phi)\big) \\
	&=\log s+{1\over 2}\log \big(1+(t/s)^2-2t/s\cos(\theta-\phi)\big)\\
	&=\log s+{1\over 2} \log\big( (1-t/s e^{i(\theta-\phi)})(1-t/s e^{-i(\theta-\phi)})   \big) \\
	&=\log s-\sum_{m=1}^\infty { t^m\over m s^m}\cos\big(m(\theta-\phi)\big).
	\end{align*}
	For the last equality, we have used the Taylor expansion of $\log(1-z)$ because $t<s$.
	Inserting it into \eqref{eq:polar} and using that
	$$\int_0^{2\pi} \int_0^{2\pi}e^{ik(\theta-\phi)}
	\cos \big(m(\theta-\phi)\big)\,\dd\theta \dd\phi
	=0 \quad\text{for}\quad m\neq k,	$$
	we see that only the terms with $m=k$ survive. Thus, \eqref{eq:polar} becomes
	\begin{align*}
	-4N\int_0^\infty\int_0^s s^{k+1}t^{k+1} {t^k\over k s^k}
	e^{-s^2-t^2} \,\dd t\dd s &=-{4N\over k }\int_0^\infty\int_0^s
	st^{2k+1} e^{-s^2-t^2}
	\,\dd t\dd s\\= -{4N\over k }\int_0^\infty\int_t^\infty
	t^{2k+1} e^{-t^2}se^{-s^2} \,\dd s\dd t &=-{2N\over k }\int_0^\infty t^{2k+1} 
	e^{-2t^2}\,\dd t
	\end{align*}
	giving the desired result by induction.
\end{proof}

\section{Preliminaries on the Bergman kernel}\label{sec-pre}

In this section, we introduce the notion of Green function and Bergman kernel associated with the line bundle $\oL^n$ and state the diagonal and off-diagonal formulas for it. Several auxiliary functions will be defined, which will be used to calculate the expectation and variance in Section \ref{sec-exp-var}.

 	There exists a Green function $G(x, y)$ of $(X, \omega)$ defined on $X \times X$, which is smooth outside the diagonal, such that for every $x,y \in X$ (c.f.\ \cite{Lang}--Theorem 2.2),
 	\begin{align}\label{Green}
 	G(x,y)=G(y,x),\qquad 
 	\int_X G(x, \cdot) \, \omega = 0 \qquad \text{and} \qquad \ddc G(x, \cdot) = \delta_x - \omega.
 	\end{align}

 The Green function has the same singularity as the logarithmic function. Let $z$ be a local normal coordinate near $x$ with $z(x)=0$. There exist a small $r_0>0$ and a smooth function $\varrho(y,w)$ on $\D(x,r_0)\times \D(x,r_0)$ such that
 \begin{equation}\label{decom-green}
 		G(y,w)=\log|z(y)-z(w)|+\varrho(y,w).
 	\end{equation} 
To simplify the notation, we will still write $y,w$ for $z(y),z(w)$ in the coordinate $z$. This should not cause any confusion. The local integrability of $\log$ function yields the following bound.

\begin{lemma}\label{lem-G-int}
	For $k\geq 1,x\in X$,  we have
	\begin{align*}
	\int_X \big|G(x,y)\big|^k \,\omega(y) \leq C_k  
	\end{align*}
	for some $C_k>0$ independent of $x$. 
\end{lemma}

Recall that the Bergman kernel of $H^0 (X,\oL^n)$ is defined by
$$ K_{N}(x,y):= \sum_{j=1}^{N} s_j(x) \otimes \overline{s_j(y)}.  $$
It is a section of $\oL^n\boxtimes \overline {\oL^n}$, with the metric induced by $\fh_n$.
Choose local unitary frames $e_n$ near $x$ and $y$. The section $e_n(x)\otimes\overline{e_n(y)}$ has norm one, and we can write
\[
K_N(x,y)=\widetilde K_N(x,y)\,e_n(x)\otimes\overline{e_n(y)},
\]
where $\widetilde K_N(x,y)$ is a scalar function in the chosen frames and $|\widetilde K_N(x,y)|$ is the pointwise norm of the kernel. For convenience, from now on, we still use $K_N(x,y)$ to denote this scalar coefficient. Since $S_n$ is assumed to be orthonormal, one has
\begin{equation}\label{kenel-density}
\rho^N(x_1,x_2,\dots,x_{N})= {1\over N !}   \det \big( K_N(x_j,x_\ell) \big)_{1\leq j,\ell\leq N} \,\omega(x_1)\omega(x_2)\cdots\omega(x_{N}).
\end{equation}

The diagonal asymptotic expansion of the Bergman kernel is well-known, see e.g.\ \cite{cat-ber,tian-jdg,zel-imrn-sze}.

\begin{lemma}\label{lem-kxx}
	As $N\to\infty$, 
	$$K_N (x,x)=N+\fb(x)+O(N^{-1})$$
	for some smooth function $\fb$ independent of $N$.
\end{lemma}

For the off-diagonal estimate, Shiffman-Zelditch
    \cite[Theorem 2.4]{shi-zel-gafa} proved that in the normal coordinate $z$ centering at $x$, for each $b>0,\ep>0$,
$$K_N\Big({z\over \sqrt N},{w\over \sqrt N} \Big)= N e^{-{\pi\over 2}  ( |z|^2+|w|^2-2z\overline w)} \Big(  1+N^{-1/2} p_{1,x}  +N^{-1}p_{2,x} +O_{b,\ep}(N^{-3/2+\ep})    \Big)  $$
for $|z|,|w|\leq b\sqrt {\log N}$, where $p_{1,x},p_{2,x}$ are polynomials in $z,\overline z,w,\overline w$ with  degree at most $10$ and  coefficients uniformly bounded in $x$. Moreover, Lu-Shiffman \cite[Theorem 2.1]{lu-shi-jga} improved the formula by showing that $p_{1,x}$ vanishes. Therefore, we have the following 
 off-diagonal formulas for the Bergman kernel.

 \begin{proposition} \label{bergman-formula}
 	Let $z$ be a normal local coordinate of $x\in X$ with respect to $\omega$ centering at $x$. Then in this coordinate, for each $b>0$ and $|z|,|w|\leq b\sqrt{\log N\over  N}$,
 	$$K_N(z,w)=  N e^{-{\pi N\over 2}  ( |z|^2+|w|^2-2z\overline w)} \Big( 1+O_b\big(N^{-1}+N^4|z|^{10}+N^4|w|^{10}\big)  \Big)$$
 	and
 	$$ |K_N(z,w)|^2=  N^2 e^{-{\pi N}   |z-w|^2} \Big( 1+O_b\big(N^{-1}+N^4|z|^{10}+N^4|w|^{10}\big)  \Big)    $$
 \end{proposition}

We  have the following bound (\cite[Theorem 2.4]{shi-zel-gafa}) for $K_N(x,y)$ when $x,y$ are far away from each other.

\begin{lemma}\label{lem-exp}
	For any $k>0$, there exist $b>0$ depending on $k$ such that if $\dist(x,y)>b\sqrt{\log N\over N}$,
	$$ |K_N(x,y)|\leq   N^{-k}. $$
\end{lemma}

We also need the following global bound for the Bergman kernel \cite[Theorem 5.7]{ber-cmp2}.

\begin{lemma}\label{lem-globalber}
	There exist positive constants $C$ and $c$ independent of $N$ such that 
	$$ |K_N(x,y)| \leq CN e^{-c\sqrt N  \dist(x,y)} \quad\text{for all}\quad x,y\in X. $$
\end{lemma}

From now on, we will fix some $b$ large, so that by Lemma \ref{lem-exp}, for the integral involving $K_N(x,y)$ below, one can cut at $\dist(x,y)\leq b\sqrt{\log N\over N}$ with an $O(N^{-2})$ error.

Recall the   reproducing property for Bergman kernel. 

\begin{lemma} \label{lem-reproducing}
	For $x_1\in X$, we have
	$$\int K_N(x_1,x_2)K_N(x_2,x_3)\cdots K_N(x_k,x_1)\,\prod_{j=2}^k\omega(x_j)=K_N(x_1,x_1).$$
In particular
	$$\int |K_N(x,w)|^2\,\omega(w)=K_N(x,x).$$
\end{lemma}

Fix a local normal coordinate $z$ centering at $x$. Assume it contains the open ball $\B(x,r_0)$ for some $0<r_0<1/2$.
We have the following error estimate based on the reproducing property.

\begin{lemma} \label{lem-logw-z}
For any point $w$ in this coordinate $z$ with $0<|w|\leq r_0$, we have
	$$ \int_{|z|\leq b\sqrt{\log N\over N}} \log|w-z|  |K_N(z,0)|^2 \,\omega(z)= N\log|w|+  {N\over 2} E_1(\pi N|w|^2) +O\big(-\log|w|\big).  $$
	Moreover, 
	$$ \int_{|z|\leq b\sqrt{\log N\over N}} \log|z|  |K_N(z,0)|^2 \,\omega(z)= -{N\log N\over 2}-{\log \pi +\gamma_E\over 2}N+O(\log N).  $$
		
\end{lemma}

\begin{proof}
	By Proposition  \ref{bergman-formula} and \eqref{defn-normal}, the left hand side is
	\begin{align}
	\int_{|z|\leq  b\sqrt{\log N\over N}}\log |w-z| N^2 e^{-\pi  N |z|^2} \Big( 1+O\big(N^{-1} + N^4 |z|^{10}\big)  \Big)
	\label{error-logz}    {i\over 2}\Big(1+O(|z|^2)\Big) \,\dd z \wedge  \dd \overline z.  
	\end{align}
	The leading term above is 
	\begin{align*}
	\int_{|z|\leq  b\sqrt{\log N\over N}}\log |w-z| N^2 &e^{-\pi N  |z|^2}   {i\over 2}  \,\dd z \wedge  \dd \overline z   
	= \int_0^{2\pi}\int_0 ^{ b\sqrt{\log N\over N}}\log|w-re^{i\theta}| N^2 e^{-\pi N r^2}   r \,\dd r\dd\theta  \nonumber\\
	&= \int_0 ^{ b\sqrt{\log N\over N}}  r  N^2 e^{-\pi N  r^2}   \int_0^{2\pi}\log|w-re^{i\theta}|  \,\dd\theta \dd r \nonumber\\
	&=2\pi\int_0 ^{ b\sqrt{\log N\over N}}  r  N^2 e^{-\pi N  r^2}    \log\big(\max(r,|w|)  \big) \,\dd r. \nonumber\\
	&\overset{N^{-2}}{\asymp}  2\pi\int_0 ^{\infty}  r  N^2 e^{-\pi N  r^2}    \log\big(\max(r,|w|)  \big) \,\dd r. \nonumber\\
	&=2\pi\int_0 ^{|w|}  r  N^2 e^{-\pi N  r^2}    \log |w|  \,\dd r+2\pi \int_{|w|} ^{\infty}  r  N^2 e^{-\pi N  r^2}    \log r \,\dd r \nonumber \\
&= N\log |w|(1-e^{-\pi N  |w|^2})  +Ne^{-\pi N  |w|^2}\log |w|+{N\over 2} E_1(\pi N|w|^2)  \nonumber\\
&=N\log|w|+{N\over 2} E_1(\pi N|w|^2),\nonumber
	\end{align*}
	where we have used Lemma \ref{lem-E1}.
	
	We now bound the error terms of \eqref{error-logz}. Notice that the integrand of the error terms has one extra $N^{-1}$, $N^4|z|^{10}$ or $|z|^2$ factor. Take $N$ large enough such that  $\log |w-z|<0$. By similar steps as above, we have
	\begin{align*}
	&\bigg|\int_{|z|\leq b\sqrt{\log N\over N}}\big(N^{-1}+N^4|z|^{10}+ |z|^2\big)\log|w-z| N^2 e^{-\pi N  |z|^2}   {i\over 2} \,\dd z \wedge  \dd \overline z\bigg|\\
	&= 2\pi\int_0 ^{ b\sqrt{\log N\over N}} -\big(N^{-1}+N^4r^{10}+ r^2\big) r  N^2 e^{-\pi N  r^2}    \log\big(\max(r,|w|)  \big) \,\dd r  \\
	&\leq -\log|w|2\pi\int_0 ^{ b\sqrt{\log N\over N}} \big(N^{-1}+N^4r^{10}+ r^2\big) r  N^2 e^{-\pi N  r^2} \,\dd r \lesssim -\log|w|.
	\end{align*}
	This proves the first assertion. The second one is similar by using the second equality in Lemma \ref{lem-E1}.
\end{proof}

\begin{lemma}\label{intGK}
	For every $x\in X$, we have uniformly on $x$,
	$$\int G(x ,y)|K_N(x,y) |^2\,\omega(y)=- {N\log N\over 2}+O(N)$$
	and
	$$\int G(x ,y)^2|K_N(x,y) |^2\,\omega(y)=O(N\log^2 N).$$
\end{lemma}

\begin{proof}
	For the first integral, by fixing a local normal coordinate centering at $x$, we have  by  Lemmas \ref{lem-exp} and \ref{lem-logw-z} that,
	\begin{align*}
	&\int G(x,y)|K_N(x,y) |^2\,\omega(y)\overset{N^{-2}}{\asymp}  \int_{\dist(x,y)\leq b\sqrt{\log N\over N}} G(x,y)|K_N(x,y) |^2\,\omega(y)  \\
	&= \int_{|z|\leq b\sqrt{\log N\over N}} \log|z| |K_N(z,0)|^2 \,\omega(z) +\int_{\dist(x,y)\leq b\sqrt{\log N\over N}} \varrho(x,y)|K_N(x,y) |^2\,\omega(y)   
	\\&=- {N\log N\over 2}+O(N).  
	\end{align*}
 The second one is similar, which is a consequence of
	$$ \int_{|z|\leq b\sqrt{\log N\over N}} \log^2|z|  N^2    e^{-\pi N  |z|^2}   {i\over 2}  \,\dd z \wedge  \dd \overline z  \lesssim N\log ^2N.    $$ 
	The proof is complete.
\end{proof}

For a local normal coordinate $z$ centering at $x$, by Taylor's formula, 
\begin{equation}
\label{diff-rho}
 \varrho( w,z)-\varrho(w,0)=c_{w,x} z+ \overline{c_{w,x} z}+ O(|z|^2).
   \end{equation}
 for some  constant $c_{w,x}\in\C$. If $w$ is also in this coordinate, then
 \begin{equation}
 \label{diff-rho2}
 \varrho( w,z)-\varrho(w,0) =   c_x z+ \overline{c_x z}+ O(|z|^2)+O(|w|^2)
 \end{equation}
 for some   constant $c_{x}\in \C$.

 For $w\neq x$ in $X$, define
 $$ F_N(w,x):=    \int G(w,y) |K_N(x,y)|^2\,\omega(y)-G(w,x) K_N(x,x).      $$
 When replacing $G$ by a smooth function, the $F_N$ is related to \textit{Toeplitz operator} and  \textit{Berezin	transform} (c.f.\ \cite{fin-cmh,fin-mann,ma-mar-toep}) and the  behavior of $F_N$ is already known.   However, here $G(x,\cdot)$ has a logarithmic singularity at $x$, where the local estimates of Bergman kernel allow us to compute the integral. Because of the singularity, the formula for Bergman kernel in a window of size $1/\sqrt N$ is no longer sufficient. One has to enlarge the window size to $\sqrt{\log N/N}$,  eliminating the singularity. The following lemma is crucial in our calculation of the variance.

\begin{lemma}\label{lem-gwy-gwx}
   When $w$ is in the coordinate $z$ and $|w|< r_0$,
   $$ F_N(w,x) = {N\over 2} E_1(\pi N|w|^2) +O\big(-\log|w|  \big).$$
   Otherwise, $F_N(w,x)=O(1)$.
\end{lemma}

\begin{proof}
We use Lemma \ref{lem-reproducing} to write
$$ F_N(w,x)=  \  \int  \big(G(w,y)- G(w,x) \big) |K_N(x,y)|^2\,\omega(y)$$
and by Lemma \ref{lem-exp}, we have
\begin{align}
    \Big|  \int_{\dist(x,y)> b\sqrt{\log N\over N}}  \big(G(w,y)- G(w,x) \big) |K_N(x,y)|^2\,\omega(y) \Big|\lesssim  N^{-2}\big(1+|G(w,x)|\big). \label{ineq1+G}
   \end{align}

\medskip
\noindent\textbf{Case 1:}  $w$ is in the coordinate $z$ and $|w|< r_0$. 
Recall from \eqref{decom-green}   $$G(w,y)-G(w,x)=\log|w-y|+\varrho(w,y)-\log|w|-\varrho(w,x)$$ in the coordinate  $z$. By \eqref{diff-rho} and 
using Proposition \ref{bergman-formula},
\begin{align*}
  &\int_{\dist(x,y)\leq   b\sqrt{\log N\over N}}  \big(\varrho(w,y)- \varrho(w,x) \big) |K_N(x,y)|^2\,\omega(y)  \\
&=
 \int_{|y|\leq b\sqrt{\log N\over N}}\big( c_{w,x} y+ \overline{c_{w,x}y}+ O(|y|^2)  \big) N^2 e^{-\pi N  |y|^2}  \Big( 1+O\big(N^{-1} + N^4 |y|^{10}\big)  \Big)\\
     &\quad\quad\qquad\qquad\qquad\qquad\qquad\qquad\qquad\qquad\qquad\qquad\,{i\over 2} \big(1+O(|y|^2)\big)\dd y \wedge  \dd \overline y.   
\end{align*}
By symmetry, the  leading term
$$\int_{|y|\leq b\sqrt{\log N\over N}}\big( c_{w,x} y+ \overline{c_{w,x}y}  \big) N^2 e^{-\pi N  |y|^2}
\,{i\over 2}\dd y \wedge  \dd \overline y=0.$$
For the error terms, we have
$$ \int_{|y|\leq b\sqrt{\log N\over N}}  |y|^2 N^2 e^{-\pi N  |y|^2}   \, {i\over 2} \dd y \wedge  \dd \overline y \overset{N^{-2}}{\asymp}  \int_0^{2\pi}\int_0 ^{\infty}  r^3 N^2 e^{-\pi N r^2} \,\dd r\dd\theta={1\over \pi} $$
and
\begin{align*} 
 &\int_{|y|\leq b\sqrt{\log N\over N}}   \big(N^{-1}+N^4|y|^{10}\big)   N^2 e^{-\pi N  |y|^2} \, {i\over 2}\dd y \wedge  \dd \overline y  \\
  &\overset{N^{-2}}{\asymp}    \int_0^{2\pi}\int_0 ^{\infty}  (N^{-1}+N^4r^{10})  N^2  e^{-\pi N r^2} r\,\dd r\dd\theta  =O(1).
\end{align*}

On the other hand, using that $\int_{|y|\leq b\sqrt{\log N\over N}}|K_N(0,y)|^2\,\omega(y)= N+O(1)$, we have by Lemma \ref{lem-logw-z}, 
\begin{align*}
\int_{|y|\leq b\sqrt{\log N\over N}} \big(\log|w-y|-\log|w| \big) |K_N(0,y)|^2\,\omega(y)
=   {N\over 2} E_1(\pi N|w|^2) +O\big(-\log|w| \big).
\end{align*}
Thus, combining this with \eqref{ineq1+G}, where $|G(w,x)|\lesssim \big|\log |w|\big|$, we conclude that in this case,
\begin{align*}
\int \big(G(w,y)- G(w,x) \big) |K_N(x,y)|^2\,\omega(y) = {N\over 2} E_1(\pi N|w|^2) +O\big(-\log|w|\big). 
\end{align*}

\medskip 
\noindent\textbf{Case 2:} Otherwise, we can find some $r_1<r_0$ such that $\dist(w,x)>r_1$. 
 Recall that the Green function $G$ is smooth away from the diagonal. The function $G(w,\cdot)$ has a bounded $\Cc^2$ norm in this case. One can replace  $\varrho$ by $G$ in the proof of Case 1, getting 
$$\int_{\dist(x,y)\leq b\sqrt{\log N\over N}}  \big(G(w,y)- G(w,x) \big) |K_N(x,y)|^2\,\omega(y) =O(1).  $$
Together with \eqref{ineq1+G}, where  $|G(w,x)|$ is bounded in this case, we finish the proof     of the lemma.
\end{proof}

\begin{lemma}\label{lem-FN}
	For any $x\in X$, we have
	$$\int \big |F_N(w,x)\big|\, \omega(w)=O(1) \quad \text{and}\quad \int F_N(w,x)^2\,\omega(w)={\log 2\over 2}N+O(\log N).$$
\end{lemma}

\begin{proof}
  According to Lemma \ref{lem-gwy-gwx},
  \begin{align*}
  &\int |F_N(w,x)|\, \omega(w)=\int_{\dist(w,x)<r_0} |F_N(w,x)|\, \omega(w)+O(1)\\
  &=\int_{|w|<r_0} \Big| {N\over 2} E_1(\pi N|w|^2) +O\big(-\log|w|  \big)\Big| \,\omega(w)+O(1)\\
  &={N\over 2} \int_{|w|<r_0}   E_1(\pi N|w|^2)\,\omega(w)+O(1) 
  \end{align*}
  and 
  \begin{align*}
  &\int F_N(w,x)^2\,\omega(w)=\int_{\dist(w,x)<r_0} F_N(w,x)^2\, \omega(w)+O(1)\\
  &= \int_{|w|<r_0} \Big( {N\over 2} E_1(\pi N|w|^2) \Big)^2 +O\big(\log^2|w|\big)+O\Big( {N\over 2}\big| \log|w| \big|E_1(\pi N|w|^2) \Big)\,\omega(w).
  \end{align*}
 Lemma \ref{lem-intE1} gives the result.  
\end{proof}

\begin{lemma}\label{lem-diffFxw}
	 We have
	$$\int F_N(w,x)F_N(x,w)\,\omega(w)\omega(x)={\log 2\over 2}N+O(\log N).  $$
\end{lemma}

\begin{proof}
	We first bound the difference between $E_1 (\pi N |w|^2)$ and $E_1 \big(\pi N \dist(w,x)^2\big)$ when $w$ is  in the normal coordinate $z$, where we have $ \dist(w,x)^2=|w|^2+O(|w|^4)$.  By mean value theorem,
	$$ N\big|E_1 (\pi N |w|^2)-E_1 \big(\pi N \dist(w,x)^2\big) \big| \lesssim N(\pi N |w|^2)^{-1} e^{-\pi N |w|^2} \cdot N|w|^4 \lesssim N|w|^2e^{-\pi N |w|^2},$$
	which is bounded.
	It follows from Lemma \ref{lem-gwy-gwx} that 
	\begin{align}\label{boundfwx}
	\big|F_N(w,x)-F_N(x,w) \big|\lesssim  \big|\log|\dist(w,x)|\big|+E_1(\pi N|\dist(w,x)|^2). 
	\end{align}
		
	Rewrite
	\begin{align*}
	F_N(w,x)F_N(x,w) ={1\over 2}  F_N(w,x)^2  +{1\over 2} F_N(x,w)^2  -{1\over 2} \big(F_N(w,x)-F_N(x,w)\big)^2.
	\end{align*}
 Lemma \ref{lem-FN} gives
	$${1\over 2}\int F_N(w,x)^2 \,\omega(w)\omega(x)+{1\over 2}\int F_N(x,w)^2 \,\omega(w)\omega(x)={\log 2\over 2}N+O(\log N).$$
	It remains to bound
	$	\int\big(F_N(w,x)-F_N(x,w)\big)^2\,\omega(w)\omega(x)$. By Lemma  \ref{lem-gwy-gwx}, we only need to consider the integral for $\dist(w,x)<r_0$, which by \eqref{boundfwx}, is
	\begin{align*}
	&\lesssim \int_{\dist(w,x)<r_0} \Big( \big|\log|\dist(w,x)|\big|+E_1(\pi N\dist(w,x)^2) \Big)^2 \,\omega(w)\omega(x)\\
	&\lesssim \int_{|w|<r_0} \Big( \big|\log|w|\big|+E_1(\pi N|w|^2) \Big)^2 \,\omega(w)\omega(x).
	\end{align*}
	This quantity is $O(1)$   by  Lemma \ref{lem-intE1}. This completes the proof of the lemma.
	\end{proof}

For $x\in X$, define
\begin{align*}
 H_N(x):=\int G(x,y)G(x,w) K_N(x,y)K_N(y,w)&K_N(w,x)\,\omega(y)\omega(w) \\
 &-  \int G^2(x,y)|K_N(x,y)|^2\,\omega(y).    
 \end{align*}

\begin{lemma} \label{lem-HN}
	We have
	$$ H_N(x)=-{\pi^2\over 24}N+O(\log^2 N).$$
\end{lemma}

\begin{proof}
	We use Lemma \ref{lem-reproducing} to rewrite 
	$$ H_N(x)=\int G(x,y)\big(G(x,w)-G(x,y)\big) K_N(x,y)K_N(y,w)K_N(w,x)\,\omega(y)\omega(w).$$
	By Lemma \ref{lem-exp}, we only need to consider the integral for $\dist(x,y),\dist(y,w)\leq b\sqrt{\log N\over N}$. By \eqref{decom-green}, \eqref{diff-rho2}, and
	applying   Proposition \ref{bergman-formula}, 
	\begin{align*}
	H_N(x)\overset{N^{-2}}{\asymp}  \int_{|y|,|w|\leq 2b\sqrt{\log N\over N}} &\big( \log|y|+\varrho(x,x)+ c_{x}y+\overline {c_xy}+O(|y|^2)\big)\\
&	\big(\log|w|-\log|y|+c_{x}  w+ \overline{c_xw}-c_xy-\overline{c_x y}+O(|w|^2)+ O(|y|^2)  \big)\\
&	N e^{-{\pi\over 2} N |y|^2} \big( 1 +O(N^{-1}+N^4|y|^{10}) \big)  \\ 
& N e^{-{\pi\over 2} N (|y|^2+|w|^2-2y\overline w)} \big( 1+ O(N^{-1}+N^4|y|^{10}+N^4|w|^{10}) \big)  \\
&	N e^{-{\pi\over 2} N |w|^2} \big( 1+O(N^{-1}+N^4|w|^{10}) \big)\\
&	\,{i\over 2}\big(1+O(|y|^2)\big) \,\dd y \wedge  \dd \overline y\,{i\over 2}\big(1+O(|w|^2)\big) \,\dd w \wedge  \dd \overline w.
	\end{align*}
	Using Lemma \ref{lem-exp} one more time, with an $O(N^{-2})$ error, we can enlarge the integral domain to $y,w\in \C$.
	After that, the leading term is
	\begin{align*}
\int \big( \log|y|+\varrho(x,x)+ c_xy+\overline {c_xy} \big) 
\big(\log|w|-\log|y|+c_{x} w+ \overline{c_xw}-c_xy-\overline {c_xy}   \big)\\
N ^3  e^{-\pi N( |y|^2+|w|^2-y\overline w)}   
\,{i\over 2} \,\dd y \wedge  \dd \overline y\,{i\over 2}  \,\dd w \wedge  \dd \overline w.
\end{align*}
After expanding the products of these two parentheses,
we have by Lemma \ref{lem-radial} and symmetry, the integrals of $$(c_{x}y+\overline{c_x y})\log|w|,\quad (c_{x} y+\overline {c_xy})\log|y|,\quad \log|y|(c_{x}  w+ \overline{c_xw}),\quad \log|y|(c_xy+\overline{c_xy})$$ are zero. By symmetry, the integrals of  $$\varrho(x,x)\big(\log|w|-\log|y | \big),\quad \varrho(x,x)(c_x w+ \overline{c_xw}-c_xy-\overline {c_xy}) $$    are also zero.   Lemma \ref{lem-ywlog} shows that the integral of
$(c_{x}y+\overline {c_xy}) (c_{x}w+ \overline{c_xw}-c_xy-\overline {c_xy})    $
is $O(1)$.
It remains to compute 
\begin{align*}
\int  \log|y|
\big(\log|w|-\log|y| \big)
N^3 e^{-{\pi} N( |y|^2+|w|^2-y\overline  w)}
\,{i\over 2} \,\dd y \wedge  \dd \overline y\,{i\over 2}  \,\dd w \wedge  \dd \overline w,
\end{align*}
which by Lemma \ref{lem-logw-logy}, is $-\pi^2N/24$.

Now we estimate the error terms. Observe that all the error integral terms are of the type in   Lemma \ref{lem-ywlog}, which says  that they are at most $O(\log^2 N)$. 
\end{proof}

For $w\in X$, define
$$ I_N(w):=\int G(x,y) K_N(x,y)K_N(y,w)K_N(w,x)\,\omega(x)\omega(y)- \int G(x,w)|K_N(x,w)|^2\,\omega(x).      $$

\begin{lemma}\label{lem-IN}
	We have
	$$I_N(w)=O(\log N).$$
\end{lemma}

\begin{proof}
	We use Lemma \ref{lem-reproducing} to rewrite 
	$$ I_N(w)=\int \big(G(x,y)-G(x,w)\big) K_N(x,y)K_N(y,w)K_N(w,x)\,\omega(x)\omega(y).$$
		By Lemma \ref{lem-exp}, we only need to restrict the integral to the region where $\dist(x,w),\dist(y,w)$ are  at most $b\sqrt{\log N\over N}$. We fix a local normal coordinate centering at $w$. By \eqref{decom-green}, \eqref{diff-rho2} and
	applying  Proposition \ref{bergman-formula}, we get 
	\begin{align*}
	I_N(w)\overset{N^{-2}}{\asymp}  \int_{|x|,|y|\leq b\sqrt{\log N\over N}} & \big(\log|x-y|-\log|x|+c_{w} y+ \overline{c_wy}+O(|x|^2)+ O(|y|^2) \big) \\  & N e^{-{\pi\over 2} N (|x|^2+|y|^2-2x\overline y)} \big( 1+ O(N^{-1}+N^4|x|^{10}+N^4|y|^{10}) \big)  \\
	&N e^{-{\pi\over 2} N |y|^2} \big( 1+ O(N^{-1}+N^4|y|^{10}) \big)\\
	&	N e^{-{\pi\over 2} N |x|^2} \big( 1 + O(N^{-1}+N^4|x|^{10})\big)\\
	&\,{i\over 2}\big(1+O(|x|^2)\big) \,\dd x \wedge  \dd \overline x \,{i\over 2}\big(1+O(|y|^2)\big) \,\dd y \wedge  \dd \overline y.
	\end{align*}
Using Lemma \ref{lem-exp} one more time, with an $O(N^{-2})$ error, we can enlarge the integral domain to $x,y\in \C$.
 For the leading term, we make the change of variable $z=x-y$ and apply  Lemma \ref{lem-radial}, getting
		\begin{align}  
		&\int  \big(\log|x-y|-\log|x|\big) N^3 e^{-\pi N (|x|^2+|y|^2-x\overline y)}  \,{i\over 2}   \,\dd x \wedge  \dd \overline x \,{i\over 2} \,\dd y \wedge  \dd \overline y  \label{logxy-logx}\\
		&=  \int  \big(\log|z|-\log|x|\big) N^3 e^{-\pi N (|x|^2+|x-z|^2-x\overline {x-z})}  \,{i\over 2}   \,\dd x \wedge  \dd \overline x \,{i\over 2} \,\dd z \wedge  \dd \overline z  \nonumber\\
		&=\int  \big(\log|z|-\log|x|\big) N^3 e^{-\pi N (|x|^2+|z|^2-z\overline {x})}  \,{i\over 2}   \,\dd x \wedge  \dd \overline x \,{i\over 2} \,\dd z \wedge  \dd \overline z   =0  . \nonumber
	\end{align}
  Moreover, by symmetry and Lemma \ref{lem-radial}, the integral of $c_{w} y+ \overline{c_wy}$ also vanishes.
 
	It remains to bound the error terms, which by Lemma \ref{lem-ywlog} after changing variable $z=x-y$ as in \eqref{logxy-logx}, they are all of $O(\log N)$. This ends the proof of the lemma.
\end{proof}

 For $x\in X$, define
	\begin{align*}  
	J_N(x):=\int G(x,y)G(w,z) K_N(x,y)K_N(y,w)K_N(w,z)K_N(z,x)\,\omega(y)\omega(w)\omega(z) \\
	-\int G(x,y)G(x,w) K_N(x,y)K_N(y,w)K_N(w,x)\,\omega(y)\omega(w). 
 \end{align*}

 \begin{lemma}\label{lem-JN}
 	We have
 	$$J_N(x)=O(\log^2 N).$$
 \end{lemma}
 
 \begin{proof}
 	We use Lemma \ref{lem-reproducing} to rewrite 
 	$$ J_N(x)=\int G(x,y)\big(G(w,z)-G(x,w)\big) K_N(x,y)K_N(y,w)K_N(w,z)K_N(z,x)\,\omega(y)\omega(w)\omega(z).$$
 	By Lemma \ref{lem-exp}, we only need to restrict the integral to the region where $\dist(x,y),\dist(x,w),\dist(x,z)$ are  at most $b\sqrt{\log N\over N}$. We fix a local normal coordinate centering at $x$. By \eqref{decom-green}, \eqref{diff-rho}, \eqref{diff-rho2} and
 	applying  Proposition \ref{bergman-formula}, we get
 	\begin{align*}
 	J_N(x)\overset{N^{-2}}{\asymp} &\int_{|y|,|w|,|z|\leq b\sqrt{\log N\over N}}  \big(\log|y|+\varrho(x,x)+c_{x}y+\overline{c_{x}y}+O(|y|^2) \big) 
 	\\
 	& \big(\log|w-z|-\log|w| +c_{x} z+ \overline{c_x z} + O(|w|^2)+O(|z|^2) \big) \\  
 	&	N e^{-{\pi\over 2} N |y|^2} \big( 1 +O(N^{-1}+N^4|y|^{10}) \big)\\
 	  &N e^{-{\pi\over 2} N (|y|^2+|w|^2-2y\overline w)} \big(1 +O(N^{-1}+N^4|y|^{10}+N^4|w|^{10}) \big)  \\
 	  &	N e^{-{\pi\over 2} N (|w|^2+|z|^2-2w\overline z)} \big( 1+O(N^{-1}+N^4|w|^{10}+N^4|z|^{10})\big)\\
 	&	N e^{-{\pi\over 2} N |z|^2} \big(1+ O(N^{-1}+N^4|z|^{10})\big)\\
 &	\,{i\over 2}\big(1+O(|y|^2)\big) \,\dd y \wedge  \dd \overline y \,{i\over 2}\big(1+O(|w|^2)\big) \,\dd w \wedge  \dd \overline w\,{i\over 2}\big(1+O(|z|^2)\big) \,\dd z \wedge  \dd \overline z.
 	\end{align*}
 	Using Lemma \ref{lem-exp} one more time, with an $O(N^{-2})$ error, we can enlarge the integral domain to $y,w,z\in \C$.
 	The leading term is 
 	\begin{align*}
 	\int  \big(\log|y|+&\varrho(x,x)+c_{x}y+\overline{c_{x}y} \big)\big(\log|w-z|-\log|w| +c_{x} z+ \overline{c_x z} \big) \\
 	& N^4 e^{-\pi N (|y|^2+|w|^2+|z|^2 -y\overline w-w\overline z)} 
 	\,{i\over 2}  \,\dd y \wedge  \dd \overline y \,{i\over 2}  \,\dd w \wedge  \dd \overline w\,{i\over 2}  \,\dd z \wedge  \dd \overline z.
 		\end{align*}
 		As in \eqref{logxy-logx}, after making the change of variable to $w-z$ and using Lemma  \ref{lem-radial}, one will see that 
   the integrals of $$ \log|y|(c_{x}z+ \overline{c_xz}),\quad \varrho(x,x)\big(\log|w-z|-\log|w|\big),\quad     \varrho(x,x)(c_xz+ \overline{c_xz}),\quad c_{x}y\big(\log|w-z|-\log|w| \big)$$ are all zero. Lemma \ref{lem-ywlog} shows that the integral of $(c_xy+\overline{c_{x}y})(c_xz+\overline{c_x z})$ is $O(1)$.  The only remaining integral in the leading term we need to compute is
 	 	\begin{align*}
 	\int\log|y| \big(\log|w-z|-\log|w|\big)   N^4 e^{-\pi N (|y|^2+|w|^2+|z|^2-y\overline w-w\overline z)}   
 	\,{i\over 2}  \,\dd y \wedge  \dd \overline y \,{i\over 2}  \,\dd w \wedge  \dd \overline w\,{i\over 2}  \,\dd z \wedge  \dd \overline z.
 	\end{align*}
   We use Lemma \ref{lem-radial} to	integrate over $y$ first, getting some constant depending on $N$ times 
 	\begin{align*} \int  \big(\log|w-z|-\log|w|\big) N^4 e^{-\pi N (|w|^2+|z|^2-w\overline z)}   
 	\,{i\over 2}  \,\dd w \wedge  \dd \overline w\,{i\over 2}  \,\dd z \wedge  \dd \overline z   =0
 	\end{align*}
 by the same reason as \eqref{logxy-logx}. We conclude that the leading term is $O(1)$.
 	
 Lastly,	
the error integral terms are  $O(\log^2 N)$ by Lemma \ref{lem-ywlog} after changing   variable   $s=w-z$. This ends the proof of the lemma.
 	\end{proof}

 For $x\in X$, define
\begin{align*}  
L_N(x):=\int G(x,y)G(w,z) K_N(x,w)K_N(w,y)K_N(y,z)K_N(z,x)\,\omega(y)\omega(w)\omega(z) \\
-\int G(x,y)G(x,w) K_N(x,y)K_N(y,w)K_N(w,x)\,\omega(y)\omega(w). 
\end{align*}

\begin{lemma}\label{lem-LN}
	We have
	$$L_N(x)=-{\pi^2\over 48}N +O(\log^2 N).$$
\end{lemma}

\begin{proof}
	We interchange $y$ and $w$   for the second term in $L_N(x)$ and
	  use Lemma \ref{lem-reproducing} to rewrite 
	$$ L_N(x)=\int G(x,y)\big(G(w,z)-G(x,w)\big) K_N(x,w)K_N(w,y)K_N(y,z)K_N(z,x)\,\omega(y)\omega(w)\omega(z).$$
	By Lemma \ref{lem-exp}, we only need to consider the integral for $\dist(x,y),\dist(x,w),\dist(x,z)\leq b\sqrt{\log N\over N}$. We fix a local normal coordinate centering at $x$. By \eqref{decom-green}, \eqref{diff-rho}, \eqref{diff-rho2} and
	applying  Proposition \ref{bergman-formula}, we get
	\begin{align*}
	L_N(x)\overset{N^{-2}}{\asymp} &\int_{|y|,|w|,|z|\leq b\sqrt{\log N\over N}}  \big(\log|y|+\varrho(x,x)+c_{x}y+\overline{c_{x}y}+O(|y|^2) \big) 
	\\
	& \big(\log|w-z|-\log|w| +c_{x} z+ \overline{c_xz}+ O(|w|^2)+O(|z|^2) \big) \\  
	&	N e^{-{\pi\over 2} N |w|^2} \big( 1+O(N^{-1}+N^4|w|^{10}) \big)\\
	&N e^{-{\pi\over 2} N (|w|^2+|y|^2-2w\overline y)} \big(1 +O(N^{-1}+N^4|w|^{10}+N^4|y|^{10}) \big)  \\
	&	N e^{-{\pi\over 2} N (|y|^2+|z|^2-2y\overline z)} \big( 1 +O(N^{-1}+N^4|y|^{10}+N^4|z|^{10})\big)\\
	&	N e^{-{\pi\over 2} N |z|^2} \big(1+O(N^{-1}+N^4|z|^{10})\big)\\
	&	\,{i\over 2}\big(1+O(|y|^2)\big) \,\dd y \wedge  \dd \overline y \,{i\over 2}\big(1+O(|w|^2)\big) \,\dd w \wedge  \dd \overline w\,{i\over 2}\big(1+O(|z|^2)\big) \,\dd z \wedge  \dd \overline z.
	\end{align*}
		Using Lemma \ref{lem-exp} one more time, with an $O(N^{-2})$ error, we can enlarge the integral domain to $y,w,z\in \C$.
	The leading term is 
	\begin{align*}
	\int  \big(\log|y|+&\varrho(x,x)+c_{x}y+\overline{c_{x}y} \big)\big(\log|w-z|-\log|w| +c_{x} z+ \overline{c_xz}\big) \\
	& N^4 e^{-\pi N (|y|^2+|w|^2+|z|^2 -w\overline y-y\overline z)} 
	\,{i\over 2}  \,\dd y \wedge  \dd \overline y \,{i\over 2}  \,\dd w \wedge  \dd \overline w\,{i\over 2}  \,\dd z \wedge  \dd \overline z.
	\end{align*}
		As in \eqref{logxy-logx}, after performing a changing of variable to $w-z$ and using Lemma  \ref{lem-radial}, one will see that the integrals of
		$$\log|y|(c_{x}z+ \overline{c_x z}),\quad \varrho(x,x)(\log|w-z|-\log|w|),\quad \varrho(x,x)(c_{x}z+ \overline{c_xz}),\quad c_{x}y  \big(\log|w-z|-\log|w|\big)$$
		are all zero. The integral of $(c_{x}y+\overline{c_{x}y}) (c_{x}z+ \overline{c_xz})$ is $O(1)$ by Lemma \ref{lem-ywlog}.
	 The only remaining integral in the leading term we need to compute is
	\begin{align*}
	\int\log|y| \big(\log|w-z|-\log|w|\big)   N^4 e^{-\pi N (|y|^2+|w|^2+|z|^2-w\overline y-y\overline z)}   
	\,{i\over 2}  \,\dd y \wedge  \dd \overline y \,{i\over 2}  \,\dd w \wedge  \dd \overline w\,{i\over 2}  \,\dd z \wedge  \dd \overline z,
	\end{align*}
	which is $-\pi^2 N/48$ by Lemma \ref{lem-logw-logwz}.

	The error terms are $O(\log^2 N)$ by Lemma \ref{lem-ywlog} after changing variable $s=w-z$. This ends the proof of the lemma.
\end{proof}

\section{Expectation and Variance}\label{sec-exp-var}

In this section, we   compute the expectation and the variance of $\Ec_N$.  For the expectation, the leading term comes from the
logarithmic singularity of $G$ at the $1/\sqrt N$ scale, where
$|K_N(x,y)|^2$ concentrates.

The second moment involves the second, third and fourth
correlation functions, according to how many distinct points
occur in the product of two Green functions. We evaluate the
corresponding integrals using the reproducing property of $K_N$,
the estimates for $F_N$, and the local integrals computed in the
preceding section. After subtracting $\big(\bE(\Ec_N)\big)^2$, the terms
larger than order $N$ cancel, leaving a variance of order $N$.

\begin{theorem}
	As $N\to\infty$,
	 $$\bE(\Ec_N)={1\over 2}N\log N+O(N).$$
\end{theorem}

\begin{proof}
	By definition and \eqref{kenel-density},
	\begin{align*}
	\bE(\Ec_N)&=\int \sum_{j\neq k} G(x_j,x_k) {1\over N !} \det \big( K_N(x_j,x_\ell) \big)_{1\leq j,\ell\leq N}  \,\omega(x_1)\omega(x_2)\cdots\omega(x_{N}) \\
	&=\int N(N-1) G(x_1,x_2) {1\over N !} \det \big( K_N(x_j,x_\ell) \big)_{1\leq j,\ell\leq N}  \,\omega(x_1)\omega(x_2)\cdots\omega(x_{N}) \\
	&= \int N(N-1) G(x_1,x_2) {(N-2)!\over N!} \det \big( K_N(x_j,x_\ell) \big)_{1\leq j,\ell \leq 2}\,\omega(x_1)\omega(x_2)\\
	&= \int G(x_1,x_2) \big( K_N(x_1,x_1)K_N(x_2,x_2)-|K_N(x_1,x_2)  |^2     \big) \,\omega(x_1)\omega(x_2).
	\end{align*}
	Here we used that $G (x,y)=G (y,x)$. The first integral, by Lemma \ref{lem-kxx}	and using that $\int_X G(x, \cdot) \, \omega = 0$, is
	\begin{align*}
	 \int G(x_1,x_2) \big( N+O(1)\big)\big( N+O(1)\big)\,\omega(x_1)\omega(x_2)
	 =O(1).
	\end{align*}
  Combining this estimate with Lemma \ref{intGK} gives the stated asymptotic formula for the expectation.
\end{proof}

To get the variance, it seems that there is no short way. One has to expand out the square of the energy sum and compute all the terms one by one.  One reason is that the higher order terms, such as $O(N^2 \log N),O(N^2),O(N\log N)$ will cancel. But they really arise in the intermediate calculations.

\begin{theorem}\label{thm-var}
	As $N\to\infty$,
	$$\Var(\Ec_N)=\Big(\log 2-{\pi^2\over 24}\Big)N+O(\log^2 N).$$
\end{theorem}

\begin{proof}
	By definition and \eqref{kenel-density},
	\begin{align*}
	\bE(\Ec_N^2)&=\int \Big(\sum_{j\neq k} G (x_j,x_k)\Big)^2 {1\over N !} \det \big( K_N(x_j,x_\ell) \big)_{1\leq j,\ell\leq N}  \,\omega(x_1)\omega(x_2)\cdots\omega(x_{N}) \\
	&= \int \Big(\sum_{j\neq k,s\neq t} G (x_j,x_k)G (x_s,x_t)\Big) {1\over N !} \det \big( K_N(x_j,x_\ell) \big)_{1\leq j,\ell\leq N}  \,\omega(x_1)\omega(x_2)\cdots\omega(x_{N}).
	\end{align*}
		The set $\{j\neq k,s\neq t\}$ can be separated into the following sets
$$
 \{ j,k,s,t \,\text{ are all distinct}\},\quad 
 \{j=s,k\neq t\},\quad
  \{  j=t, k\neq s\},$$
   $$ 
	 \{j\neq s,k= t\},\quad
	\{  j\neq t, k= s\}, \quad  \{j=s, k=t\}, \quad
	  \{j=t,k=s\}.  
$$
The number of pairs of the first category is ${N!\over (N-4)!}$. The next four categories contain ${N!\over (N-3)!}$ pairs each. The last two contain $N(N-1)$ pairs each.
According to such classification, and using  that $G (x,y)=G (y,x)$, we obtain
	$$\bE(\Ec_N^2)=2A_2+4A_3+A_4, $$
	where
	\begin{align*} 
	A_2&=\int   G (x_1,x_2) G (x_1,x_2)  \det \big( K_N(x_j,x_\ell) \big)_{1\leq j,\ell\leq 2}  \,\omega(x_1)\omega(x_2)\\
	A_3 &=\int   G (x_1,x_2)G (x_1,x_3)   \det \big( K_N(x_j,x_\ell) \big)_{1\leq j,\ell\leq 3}  \,\omega(x_1)\omega(x_2)\omega(x_3) \\
	A_4&= \int  G (x_1,x_2)G (x_3,x_4)  
	  \det \big( K_N(x_j,x_\ell) \big)_{1\leq j,\ell\leq 4}  \,\omega(x_1)\omega(x_2)\omega(x_3)\omega(x_4).
	\end{align*}
	By Propositions \ref{prop-A2}, \ref{prop-A3}, \ref{prop-A4} below, all terms cancel with each other except $N$ and $O(\log^2 N)$ terms.
	\end{proof}

The remaining part of this section is devoted to evaluate all the $A_2,A_3,A_4$.
	
		\begin{proposition}\label{prop-A2}
		As $N\to \infty$, 
		\begin{align*} 
		A_2= \int G (x,y)^2   K_N(x,x)K_N(y,y)\,\omega(x)\omega(y)-\int  G (x,y)^2|K_N(x,y)  |^2     \,\omega(x)\omega(y)
		\end{align*}
	\end{proposition}
	
	\begin{proof}
		Use definition.
	\end{proof}

	\begin{proposition}\label{prop-A3}
	 As $N\to \infty$,	
	 \begin{align*}
	 A_3&= \int  \left( \int G (x,y)\big(K_N(y,y)-N\big)\, \omega(y)  \right)^2  K_N(x,x)\,\omega(x) \\
	 &\qquad -2\int G (x,y)|K_N(x,y)|^2 G (x,w)\fb(w)  \,\omega(x)     \omega(y)  \omega(w)\\
	 &\qquad  -\int G (x,y)^2 K_N(x,x)K_{N}(y,y) \,\omega(x)\omega(y)- \int G (x,y) F_N(x,y) K_N(x,x)\,\omega(x)\omega(y)\\
	 &\qquad  +2\int G (x,y)^2 |K_N(x,y)|^2 \,\omega(x)\omega(y)-{\pi^2\over 12}N+O(\log^2 N)
	 \end{align*}
	\end{proposition}

    \begin{proof} 
	Expanding the determinant  $\det \big( K_N(x_j,x_\ell) \big)_{1\leq j,\ell\leq 3}$  gives
	$$K_{11}K_{22}K_{33}+K_{12}K_{23}K_{31}+K_{13}K_{21}K_{32}-K_{13}K_{22}K_{31}-K_{12}K_{21}K_{33}-K_{11}K_{23}K_{32}, $$
	where we abbreviate $K_N(x_j,x_\ell)$ as $K_{j\ell}$.  The proposition will follow by   Lemmas \ref{lem-a3-1}--\ref{lem-3.15} below.
	\end{proof}

   \begin{lemma} \label{lem-a3-1}
   	As $N\to\infty$,
   	\begin{align*}
   	  \int   G(x_1,x_2)G (x_1,x_3) & K_{11}K_{22}K_{33}  \,\omega(x_1)\omega(x_2)\omega(x_3)       \\&=\int  \left( \int G (x,y)\big(K_N(y,y)-N\big)\, \omega(y)  \right)^2  K_N(x,x)\,\omega(x).
   	  \end{align*}    
   	\end{lemma}

    \begin{proof}
By Lemma  \ref{lem-kxx} and using that $\int_X G(x,\cdot)\,\omega=0$,
	the left hand side is 
		\begin{align*}
	&   \int   G (x_1,x_2)G (x_1,x_3)  K_{11}\big(N+\fb(x_2)+O(N^{-1})\big) \big(N+\fb(x_3)+O(N^{-1})\big) \,\omega(x_1)\omega(x_2)\omega(x_3)     \\     
	&=\int   G (x_1,x_2)G (x_1,x_3)  K_{11}\big(K_N(x_2,x_2)-N\big) \big(K_N(x_3,x_3)-N\big) \,\omega(x_1)\omega(x_2)\omega(x_3),
	  \end{align*}
	  giving the lemma.
 \end{proof}

	 \begin{lemma}
	 	 	As $N\to\infty$,
	  	\begin{align*}   
	  	&\int   G (x_1,x_2)G (x_1,x_3)  K_{12}K_{21}K_{33}  \,\omega(x_1)\omega(x_2)\omega(x_3)      \\
	  	&= \int   G (x_1,x_2)G (x_1,x_3)  K_{13}K_{22}K_{31}  \,\omega(x_1)\omega(x_2)\omega(x_3)    \\  
	  	&=\int G (x,y)|K_N(x,y)|^2 G (x,w)\fb(w)  \,\omega(x)     \omega(y)  \omega(w)+O(\log N). \end{align*}
	  \end{lemma}
  
  \begin{proof}The first equality is by symmetry.
	  By Lemmas  \ref{lem-kxx}, \ref{lem-logw-z} and using that $\int_X G(x,\cdot)\,\omega=0$, the integral is 
	  	\begin{align*}
	  	&   \int   G (x_1,x_2)  |K_N(x_1,x_2)|^2\,\omega(x_2)  \int G (x_1,x_3)  \big(N+\fb(x_3)+O_{x_3}(N^{-1})\big) \,\omega(x_3)     \omega(x_1)   \\
	  	& =  \int   G (x_1,x_2)  |K_N(x_1,x_2)|^2\,\omega(x_2)  \int G (x_1,x_3)  \big(\fb(x_3)+O_{x_3}(N^{-1})\big) \,\omega(x_3)     \omega(x_1)   \\
	  	&= \int   G (x_1,x_2)  |K_N(x_1,x_2)|^2  G (x_1,x_3)  \fb(x_3)\,\omega(x_2)\omega(x_3)     \omega(x_1)  \\
	  	&\qquad + \int   G (x_1,x_2)G (x_1,x_3)   |K_N(x_1,x_2)|^2   O_{x_3}(N^{-1})    \,\omega(x_1) \omega(x_2)\omega(x_3).
	  	\end{align*}
	 It remains to estimate the second integral. We use Lemma \ref{intGK} to integrate $x_2$ first, and get the desired bound $O(\log N)$ since $G(x_1,x_3)$ is integrable over $\omega\otimes \omega$.
	  \end{proof}

	  \begin{lemma} \label{lem-k11k23k32}
	  	As $N\to\infty$,
	  	\begin{align*}
	  	\int   G (x_1,&x_2) G (x_1,x_3)  K_{11}K_{23}K_{32}\,\omega(x_1)\omega(x_2)\omega(x_3)=\\
	  	& \int G (x,y)^2 K_N(x,x)K_{N}(y,y) \,\omega(x)\omega(y)+ \int G (x,y) F_N(x,y) K_N(x,x)\,\omega(x)\omega(y). 
	  		\end{align*}
	  \end{lemma}
	  
	  \begin{proof}
	  	By definition of $F_N$.
  \end{proof}

  	\begin{lemma}\label{lem-3.15}
  			As $N\to\infty$,
  \begin{align*}
  &\int  G (x_1,x_2)G (x_1,x_3)  K_{12}K_{23}K_{31}\,\omega(x_1)\omega(x_2)\omega(x_3)\\
  &=\int  G (x_1,x_2)G (x_1,x_3)  K_{13}K_{32}K_{21}\,\omega(x_1)\omega(x_2)\omega(x_3)\\
  &=\int G (x,y)^2 |K_N(x,y)|^2 \,\omega(x)\omega(y)-{\pi^2\over 24}N+O(\log^2 N).
  \end{align*}
  	\end{lemma}
  
  \begin{proof}
 The first equality is by symmetry. The second one is  by Lemma \ref{lem-HN}.
\end{proof}

\begin{proposition} \label{prop-A4}
	As $N\to\infty$, 
	\begin{align*} 
	A_4&= \big(\bE(\Ec_N) \big)^2-4\int  \left( \int G (x,y)\big(|K_N(y,y)|-N\big)\, \omega(y)  \right)^2  K_N(x,x)\,\omega(x)\\
	&\qquad +2\int G (x,y)^2 K_N(x,x)K_N(y,y) \,\omega(x)\omega(y)\\&\qquad+4\int G(x,y) F_N(x,y) K_N(x,x) \,\omega(x)\omega(y)+\log 2 \cdot N\\
	&\qquad +8\int G(x,y)|K_N(x,y)|^2 G(y,w) \fb(w) \,\omega(x )\omega(y)\omega(w)\\
	&\qquad -6\int G (x,y)^2 |K_N(x,y)|^2\,\omega(x )\omega(y) +{7\pi^2\over 24}N+O(\log^2 N).
	\end{align*}
\end{proposition}

\begin{proof}
	For $\sigma\in\Sym(4)$, define
	$$ A_\sigma:=  \int   G(x_1,x_2)G (x_3,x_4)  
	\prod_{j=1}^4 K_N(x_j,x_{\sigma(j)})  \,\omega(x_1)\omega(x_2)\omega(x_3)\omega(x_4).   $$
	Then by definition,
	$$ A_4=  
	\sum_{\sigma\in\Sym(4)}\sgn(\sigma)A_\sigma.  $$

	Notice that
	\begin{align*}
	&\big(\bE(\Ec_N) \big)^2 = \bigg(\int G (x_1,x_2) \big( K_N(x_1,x_1)K_N(x_2,x_2)-|K_N(x_1,x_2)  |^2     \big) \,\omega(x_1)\omega(x_2)\bigg)^2 \\
	&=\int G (x_1,x_2) \big(K_{11}K_{22}-|K_{12}|^2\big) \,\omega(x_1)\omega(x_2)\int G (x_3,x_4) \big(K_{33}K_{44}-|K_{34}|^2\big) \,\omega(x_3)\omega(x_4) \\
	&=A_{\id}-A_{(3,4)}-A_{(1,2)}+A_{(1,2)(3,4)}.
	\end{align*}
	Thus, we do not need to calculate these four terms, as they will cancel with $\big(\bE(\Ec_N) \big)^2$ when calculating the variance. The proposition will follow by Lemmas \ref{lem-A4-1}--\ref{lem-A4-5} below.
\end{proof}

\begin{lemma}\label{lem-A4-1}
	As $N\to\infty$, 
	\begin{align*}
	&A_{(1,3)}=A_{(1,4)}=A_{(2,3)}=A_{(2,4)}\\
	&=\int  \left( \int G (x,y)\big(|K_N(y,y)|-N\big)\, \omega(y)  \right)^2  K_N(x,x)\,\omega(x)+O(1). 
	\end{align*}
\end{lemma}

\begin{proof}
	We only need to calculate $A_{(1,3)}$ by symmetry.
 By definition,
 $$A_{(1,3)} =\int   G (x_1,x_2)G (x_3,x_4)  
 |K_{13}|^2K_{22} K_{44}  \,\omega(x_1)\omega(x_2)\omega(x_3)\omega(x_4).$$
 For every fixed $x_1,x_3$, by Lemma \ref{lem-kxx} and using that $\int_X G(x,\cdot)\,\omega=0$, the integrand above is
\begin{align*}
 |K_{13}|^2\int   G (x_1,x_2)  
 \big(K_N(x_2,x_2)-N\big) \,\omega(x_2)   \int G (x_3,x_4)\big(K_N(x_4,x_4)-N\big)    \,\omega(x_4).
\end{align*}
 We use $F_N$ to write
\begin{align*}
&\int     G (x_1,x_2)  |K_{13}|^2
   \,   \omega(x_1) 
 =  G (x_2,x_3) K_{33}  +F_N(x_2,x_3). 
\end{align*}
It follows that $A_{(1,3)}$ is 
\begin{align*}
&\int G (x_2,x_3) K_{33} \big(K_N(x_2,x_2)-N\big)  G (x_3,x_4)\big(K_N(x_4,x_4)-N\big)  \,\omega(x_2)\omega(x_3)\omega(x_4)\\
& \quad +\int F_N(x_2,x_3)  \big(K_N(x_2,x_2)-N\big)G (x_3,x_4)\big(K_N(x_4,x_4)-N\big)   \,\omega(x_2)\omega(x_3)\omega(x_4).
\end{align*}
The first integral is the main term we want and  the second integral is $O(1)$ by Lemma \ref{lem-FN}.
\end{proof}

\begin{lemma}
	As $N\to\infty$, 
	\begin{align*}
	A_{(1,3)(2,4)}=A_{(1,4)(2,3)}&=\int G (x,y)^2 K_N(x,x)K_N(y,y) \,\omega(x)\omega(y)\\
	&\quad +2\int G(x,y) F_N(x,y) K_N(x,x) \,\omega(x)\omega(y)+{\log 2\over 2}N+O(\log N).\end{align*}
\end{lemma}

\begin{proof}
	By symmetry, we only need to calculate $A_{(1,4)(2,3)}$, which
by definition, is
$$ \int   G (x_1,x_2)G (x_3,x_4)  
|K_{14}|^2 |K_{23}|^2 \,\omega(x_1)\omega(x_2)\omega(x_3)\omega(x_4).$$
We use $F_N$ to write
 $$\int G (x_1,x_2)  
 |K_{23}|^2\,\omega(x_2)= G (x_1,x_3) K_{33}+F_N(x_1,x_3) $$
 and
 $$\int G (x_3,x_4)  
 |K_{14}|^2\,\omega(x_4)= G (x_3,x_1) K_{11}+F_N(x_3,x_1).$$
Therefore, $A_{(1,4)(2,3)}$ is
\begin{align*}
&\int G (x_1,x_3)^2 K_{11}K_{33} \,\omega(x_1)\omega(x_3)+ \int G(x_1,x_3) F_N(x_1,x_3) K_N(x_1,x_1) \,\omega(x_1)\omega(x_3)+ \\
& \int G(x_1,x_3) F_N(x_3,x_1) K_N(x_3,x_3) \,\omega(x_1)\omega(x_3) +\int F_N(x_1,x_3)F_N(x_3,x_1) \,\omega(x_1)\omega(x_3).
\end{align*}
From the definition of $F_N$ and using that $G(x_1,x_3)=G(x_3,x_1)$, we see that
the second and third integrals are the same. The lemma will follow by Lemma \ref{lem-diffFxw}. 
\end{proof}

\begin{lemma}
	If $\sigma$ contains a $3$-cycle, then as $N\to\infty$,
  $$A_{\sigma}=\int G(x,y)|K_N(x,y)|^2 G(y,w) \fb(w) \,\omega(x )\omega(y)\omega(w)+O(\log N).$$ 
\end{lemma}

\begin{proof}
	We only give the proof for $\sigma=(1,2,3)$. The others can be treated in the same way.
	By definition and Lemma \ref{lem-kxx}, $A_{(1,2,3)}$ is
	\begin{align*}
	&=\int   G (x_1,x_2)G (x_3,x_4)  
	K_{12}K_{23}K_{31}K_{44}  \,\omega(x_1)\omega(x_2)\omega(x_3)\omega(x_4) \\
	&= \int   G (x_1,x_2)K_{12}K_{23}K_{31} \Big( \int G (x_3,x_4) K_{44}\, \omega(x_4) \Big)
	  \,\omega(x_1)\omega(x_2)\omega(x_3) \\
	  &= \int   G (x_1,x_2)K_{12}K_{23}K_{31} \Big( \int G (x_3,x_4) \big(\fb(x_4)+O_{x_4}(N^{-1})\big)\, \omega(x_4) \Big)
	  \,\omega(x_1)\omega(x_2)\omega(x_3)\\
	  &=\int   G (x_1,x_2)K_{12}K_{23}K_{31} G (x_3,w)\fb(w)\omega(w)\omega(x_1)\omega(x_2)\omega(x_3) \\
	 &\qquad\qquad\qquad +\int   G (x_1,x_2)K_{12}K_{23}K_{31} \,O_{x_3}(N^{-1})\,\omega(x_1)\omega(x_2)\omega(x_3).
	\end{align*}

By Lemma \ref{lem-IN}, for each $x_3$,
$$\int   G (x_1,x_2) K_{12}K_{23}K_{31} \,\omega(x_1)\omega(x_2)= \int G(x_1,x_3)  |K_{13}|^2 \omega(x_1)+O(\log N).$$
   We conclude from  Lemma \ref{intGK} that  
   	\begin{align*}
   A_{(1,2,3)}&=\int G(x_1,x_3)|K_{13}|^2 G(x_3,x_4) \fb(x_4) \,\omega(x_1)\omega(x_3)\omega(x_4)  +O(\log N),
   	\end{align*}
   	giving the lemma.
\end{proof}

\begin{lemma}
	 As $N\to\infty$,
	 \begin{align*}
	 A_{(1,2,3,4)}=A_{(1,2,4,3)}&=A_{(1,3,4,2)}=A_{(1,4,3,2)}  \\
	 &=\int G (x,y)^2 |K_N(x,y)|^2 \,\omega(x)\omega(y)-{\pi^2\over 24}N+O(\log^2 N).
	 \end{align*}
\end{lemma}

 \begin{proof}
 Using Lemmas \ref{lem-HN} and \ref{lem-JN}.
\end{proof}

\begin{lemma}\label{lem-A4-5} 
	As $N\to\infty$,
	\begin{align*}
	A_{(1,3,2,4)}=A_{(1,4,2,3)}=\int G (x,y)^2 |K_N(x,y)|^2 \,\omega(x)\omega(y)-{\pi^2\over 16}N+O(\log^2 N).
	\end{align*}
\end{lemma}

\begin{proof}
	Using Lemmas \ref{lem-HN} and \ref{lem-LN}.
\end{proof}

\section{Moments and Cumulants}

In this section, we  recall the notion of joint cumulants in probability theory, and derive a formula for the cumulants of the Green energy. The reader may refer to  the surveys \cite{han-book,cum-survey} for the details.

By a \textit{partition} of a finite set $Q$, we mean a set $\Qc$, whose elements are disjoint non-empty subsets of $Q$ such that their union  is $Q$.
Let $\Qc$ be a partition of $Q$. For $s\in Q$, let $\Qc(s)$ be the element in $\Qc$ which contains $s$.

Let $Y_1,Y_2,\dots,Y_r$  be real random variables whose products  are integrable. Their \textit{joint cumulant} is defined by
\begin{align}\label{formula-cum}
 \kappa( Y_1,Y_2,\dots,Y_r):= \sum_{\Qc\in \oP([r]) } (-1)^{|\Qc|-1} (|\Qc|-1)! \prod_{I\in \Qc} \bE\Big( \prod_{s\in I} Y_s \Big).
\end{align}  
Here, $[r]:=\{1,2,\dots,r\}$, $\oP([r])$ means all the partitions of the set $[r]$ and $|\Qc|$ denotes the number of elements in $\Qc$.  For a single random variable $Y$, we define its $r$-th cumulant by 
$$ \kappa_r(Y):=\kappa(Y,Y,\dots,Y)   \quad\text{with $r$ copies},  $$
which gives
\begin{align}\label{formula-cum-2}
\kappa_r(Y):= \sum_{\Qc\in \oP([r]) } (-1)^{|\Qc|-1} (|\Qc|-1)!   \prod_{I\in \Qc} \bE ( Y^{|I|} ).
\end{align}

The inverse formula is given by
\begin{align*}
\bE \Big(\prod_{j=1}^r Y_j \Big)=\sum_{\Qc\in \oP([r]) }  \prod_{I\in \Qc} \kappa( Y_s:\,s\in I ),
\end{align*}
where $\kappa( Y_s:\,s\in I )$ means the joint cumulant of  $Y_s$ for $s\in I$, which is well-defined since joint cumulants are symmetry in their arguments.

For two partitions $\Qc_1$ and $\Qc_2$ of $Q$, we define an ordering as follows: $\Qc_1\leq \Qc_2$  if for any $A\in \Qc_1$, one has $A\subset B$ for some $B\in \Qc_2$. This means $\Qc_1$ is a refinement of $\Qc_2$. Denote by $\Qc_1 \vee \Qc_2$   the minimal partition of $Q$ such that $\Qc_1\leq \Qc_1 \vee \Qc_2$ and $\Qc_2\leq \Qc_1 \vee \Qc_2$. For example, let $\Qc_1=\big\{\{1,2\},\{3\},\{4,5\},\{6\}  \big\}$ and $\Qc_2=\big\{\{1\},\{2,3\},\{4,5\},\{6\}  \big\}$. Then $\Qc_1\vee \Qc_2= \big\{\{1,2,3\},\{4,5\},\{6\}  \big\}$.

Suppose $\Qc$ and $\Lc$ are  partitions of disjoint sets $Q$ and $L$ respectively. Then $\Qc \cup \Lc:=\{A:\,A\in\Qc \,\text{ or } \, A\in\Lc \}$ will be a partition of $Q\cup L$. If $F\subset Q$, $\Qc|_F:=\{A\cap F:\, A\in \Qc,A\cap F\neq\varnothing\}$ is a partition of $F$.

For each $\sigma\in \Sym(\Qc)$, denote by $\Qc\otimes\sigma$ the new partition obtained by merging the sets   $A\in\Qc$ belonging to the same cycle of $\sigma$. In other words, $\Qc\otimes\sigma$ is the minimal partition, where each element is invariant under $\sigma$.     For instance, if $\Qc=\big\{\{1,2\},\{3\},\{4,5\},\{6\}  \big\}$, which is a partition of $[6]$. Let $\sigma$ interchanges $\{1,2\}$ and $\{3\}$ and fixes $\{4,5\}$ and $\{6\}$, then $\Qc\otimes\sigma=\big\{\{1,2,3\},\{4,5\},\{6\}  \big\}$.

Fix the partition $$ \Rc_k:=\big\{ \{1,2\},\{3,4\},\dots,\{2k-1,2k \}   \big\}$$ of $[2k]$ and the trivial partition $\{[2k] \}$ of $[2k]$.

 As $G$ is not defined on the diagonal, we assign its value to be zero by defining a new function
$$\tG(x,x):=0  \quad\text{and} \quad  \tG(x,y):=G(x,y)  \quad\text{for} \quad x\neq y\in X.    $$
Put
$$ \widetilde\Ec_N:=\sum_{j\neq k}\tG(x_j,x_k)=\sum_{j,k}\tG(x_j,x_k). $$
 We have the following two results  for the moments and cumulants of the Green energy, see also  \cite[Proposition 3.1]{cum-survey}.

\begin{lemma}\label{lem-moment}
		For every $k\geq 1,N\geq 2k$, we have
		$$ \bE(\Ec_N^k)= \sum_{\Qc \in \oP([2k])} \sum_{\sigma\in\Sym(\Qc)} \sgn(\sigma) \int_{X^{\Qc}}   \prod_{j=1}^k\tG\big(x_{\Qc(2j-1)},x_{\Qc(2j)}\big)      \prod_{A\in\Qc} K_N(x_A,x_{\sigma(A)}) \, \prod_{A\in\Qc}\omega(x_A).$$
\end{lemma}

\begin{proof}
		By definition, $\bE(\Ec_N^k)=\bE(\widetilde\Ec_N^k)$ for every $k$, where we can write
	$$ \widetilde\Ec_N^k=\sum_{\eta:[2k]\to[N]}  \prod_{j=1}^k\tG\big(x_{\eta(2j-1)},x_{\eta(2j)}\big)  $$
	and  using \eqref{kenel-density}, we have
	\begin{align} \label{etg-1}
	\bE(\widetilde\Ec_N^k)=\sum_{\eta:[2k]\to[N]} \int_{X^{N}}   \prod_{j=1}^k\tG\big(x_{\eta(2j-1)},x_{\eta(2j)}\big)   {1\over N !}   \det \big( K_N(x_j,x_\ell) \big)_{1\leq j,\ell\leq N} \,\prod_{\ell=1}^{N}\omega(x_\ell),
	\end{align}
	where $\eta$ runs over all the maps from $[2k]$ to $[N]$.   
	Every $\eta$ will induce a partition $\Qc_\eta$  of $[2k]$ according to the value  $\eta(s)$ for $s\in[2k]$, where we put together the $s$ who has the same $\eta$-value. Note that $\Qc_{\eta_1}$ may equal  $\Qc_{\eta_2}$ for $\eta_1\neq \eta_2$,  and the number of $\eta$'s given the same partition $\Qc_\eta$ is $N (N-1)\cdots(N-|\Qc_\eta|+1)$. Moreover, the same partition will contribute the same integral in the summation \eqref{etg-1}. It follows that $\bE(\widetilde\Ec_N^k)$ is
	\begin{align*}
	&\sum_{\Qc \in \oP([2k])}{N!\over (N-|\Qc|)!} \int_{X^{N}}   \prod_{j=1}^k\tG\big(x_{\Qc(2j-1)},x_{\Qc(2j)}\big)   {1\over N !}   \det \big( K_N(x_j,x_\ell) \big)_{1\leq j,\ell\leq N} \,\prod_{\ell=1}^{N}\omega(x_\ell)\\
	&=\sum_{\Qc \in \oP([2k])} {N!\over (N-|\Qc|)!} \int_{X^{|\Qc|}}   \prod_{j=1}^k\tG\big(x_{\Qc(2j-1)},x_{\Qc(2j)}\big) \\ &\qquad\qquad\qquad\qquad\qquad\qquad\qquad\qquad { (N-|\Qc|)! \over N!}   \det \big( K_N(x_j,x_\ell) \big)_{1\leq j,\ell\leq |\Qc|} \,\prod_{\ell=1}^{|\Qc|}\omega(x_\ell)\\
	&=\sum_{\Qc \in \oP([2k])}  \int_{X^{|\Qc|}}   \prod_{j=1}^k\tG\big(x_{\Qc(2j-1)},x_{\Qc(2j)}\big)     \det \big( K_N(x_j,x_\ell) \big)_{1\leq j,\ell\leq |\Qc|} \,\prod_{\ell=1}^{|\Qc|}\omega(x_\ell).
	\end{align*}
	The result follows after expanding the determinant.
\end{proof}

For   $I\subset [k]$, put $2I:=\cup_{j\in I}\{ 2j-1,2j\}$. This $2I$ does not mean doubling the elements in $I$ in the usual sense. For every partition $\Qc$ of $2I$ and $\sigma\in \Sym(\Qc)$, define the integral
$$ \oJ(\Qc,\sigma):= \sgn(\sigma)  \int_{X^{\Qc}} \prod_{j\in I} \tG\big(x_{\Qc(2j-1)},x_{\Qc(2j)}\big)  \prod_{A\in\Qc} K_N(x_A ,x_{\sigma(A)})\,\prod_{A\in\Qc}\omega(x_A).    $$

\begin{proposition}\label{prop-cum}
	For every $k\geq 2,N\geq 2k$, we have
  \begin{align*}
  \kappa_k(\Ec_N)=\sum_{\Qc\in \oP([2k])}\sum_{\sigma\in \Sym(\Qc)\atop(\Qc\otimes \sigma)\vee \Rc_k=\{[2k]\}} \oJ(\Qc,\sigma).
  \end{align*}	
\end{proposition}

\begin{proof}
  By Lemma \ref{lem-moment}, for every $I\subset [k]$,
  $$  \bE(\Ec_N^{|I|}) = \sum_{\Lc \in \oP(2I)} \sum_{\tau\in\Sym(\Lc)}
   \oJ(\Lc,\tau).$$
	Hence  \eqref{formula-cum-2} gives
	\begin{align} \label{kec-1}
	\kappa_k(\Ec_N)=	\sum_{\Gc\in \oP([k]) } (-1)^{|\Gc|-1} \big(|\Gc|-1\big)!\prod_{I\in \Gc} \bigg(\sum_{\Lc \in \oP(2I)} \sum_{\tau\in\Sym(\Lc)}
	\oJ(\Lc,\tau) \bigg) 
	\end{align}
	
	Fix a  $\Gc\in \oP([k])$. For each $\oI_j\in\Gc$,  taking   $\Lc_j\in \oP(2I_j),\tau_j\in \Sym(\Lc_j) $ and gluing them, will give a partition  $\Qc=\cup \Lc_j$ of $[2k]$ and a permutation $\sigma=\prod \tau_j\in \Sym(\Qc)$. Note that $\sgn(\sigma)=\prod\sgn(\tau_j)$. By Fubini's theorem, we have
	\begin{align} \label{kec-2}
	\oJ(\Qc,\sigma)=\prod \oJ(\Lc_j,\tau_j). 
	\end{align}
	
 It remains to count how many different $\Gc$'s will produce the same $(\Qc,\sigma)$ in the way above. More precisely, for a  fixed pair $(\Qc,\sigma)$, where $\Qc$ is a partition of $[2k]$ and $\sigma\in \Sym(\Qc)$, we need to find all $\Gc\in\oP([k])$ such that there exists $\Lc_j\in \oP(2I_j)$ and $\tau_j\in\Sym(\Lc_j)$ for every $I_j\in\Gc$ such that 
 \begin{align}\label{produce}
 \Qc=\bigcup \Lc_j  \quad\text{and}\quad \sigma=\prod \tau_j\in \Sym(\Qc).
 \end{align} 
 It is straightforward to check  that a necessary condition is that for every $A\in\Qc$, it should be contained in  $2I$ for some $I\in \Gc$, and  $\sigma(A)$ should be contained in this $2I$ also. This condition is also sufficient because $[2k]$ is the disjoint union of $A\in\Qc$ and $2I_j$'s respectively. 
 
 \medskip
 \noindent
 \textbf{Claim:} $\Gc$ produces $(\Qc,\sigma)$ if and only if for any $B\in (\Qc\otimes\sigma)\vee\Rc_k$, $B\subset 2I$ for some $I\in\Gc$.

\begin{proof}[Proof of Claim]
	Suppose $\Gc$ produces $(\Qc,\sigma)$, then we can find $\Lc_j\in \oP(2 I_j),\tau_j\in \Sym(\Lc_j)$ for every $I_j\in \Gc$ such that \eqref{produce} holds.  To get the partition $(\Qc\otimes\sigma)\vee\Rc_k$, one needs to do two steps. The first step is to merge the elements of $\cup \Lc_j$ according to the actions of $\sigma$. After that, we merge the elements further to make $2j-1$ and $2j$ into one single element. Both of these two steps do not move elements in $[2k]\setminus 2I$ into the partition of $2I$ for any $I\in\Gc$. For $B\in (\Qc\otimes\sigma)\vee\Rc_k$, all the elements merging to $B$ should be contained in the same $2I_j$ for some $I_j\in \Gc$, which means $B\subset 2I$ for some $I\in\Gc$.
	
	Conversely, suppose for any $B\in (\Qc\otimes\sigma)\vee\Rc_k$, $B\subset 2I$ for some $I\in\Gc$. For each $I\in\Gc$, define 
	$ \Lc_I:=\{ A\in \Qc:\, A\subset 2I\}$.  The assumption implies that $\sigma(A)\in \Lc_I$ for all $A\in \Lc_I$, and hence $\tau_I:=\sigma|_{\Lc_I}$ is well-defined. One can check $(\Lc_I,\tau_I)$ satisfies condition \eqref{produce}. This ends the proof of the Claim.
\end{proof}

By the Claim, we see that if  $\Gc$ produces  $(\Qc,\sigma)$, then $(\Qc\otimes\sigma)\vee\Rc_k \leq  2\Gc$, where $2\Gc$ is the partition of $[2k]$ obtained by $2I_j$'s for $I_j\in \Gc$.  
 Therefore,  using that $|\Gc|=|2\Gc|$, we conclude from \eqref{kec-1} and \eqref{kec-2} that 
$$\kappa_k(\Ec_N)= \sum_{\Qc\in \oP([2k])}\sum_{\sigma\in \Sym(\Qc) }\sum_{\Fc\in\oP([2k])\atop  (\Qc\otimes\sigma)\vee\Rc_k \leq  \Fc}  (-1)^{|\Fc|-1} (|\Fc|-1)! \oJ(\Qc,\sigma).  $$

Let $n_{\Qc,\sigma}$ be the number of elements in $(\Qc\otimes\sigma)\vee\Rc_k$. The number of choices of $\Fc$ is exactly $|\oP([n_{\Qc,\sigma}])|$. Thus,
$$ \sum_{\Fc\in\oP([2k])\atop  (\Qc\otimes\sigma)\vee\Rc_k \leq \Fc}  (-1)^{|\Fc|-1} (|\Fc|-1)!= \sum_{ \Tc \in \oP([n_{\Qc,\sigma}])  } (-1)^{|\Tc|-1} (|\Tc|-1)!,$$
whose value is $1$ when $n_{\Qc,\sigma}=1$ and zero otherwise (\cite[(3.12)]{cum-survey}). So the only term surviving is when $(\Qc\otimes\sigma)\vee\Rc_k=\{[2k]\}$. This completes the proof of the proposition.
\end{proof}

Let $f$ be a complex-valued  function on  $X^r$ such that 
$$  \int_{X^\Lc} \big|  f(x_{\Lc(1)},x_{\Lc(2)},\dots,x_{\Lc(r)}  )\big|\, \prod_{A\in\Lc}\omega(x_A)<+\infty  $$
for any $\Lc\in\oP([r])$. Denote by $\Cyc(\Lc)$  the set of permutations in $\Sym(\Lc)$ with a single cycle of length $|\Lc|$. We have the following result of cumulant reduction analogous to \cite[Proposition 4.12]{man-mom}, see also  \cite[(2.7)]{sos-com} and \cite[Lemma 1]{feng-goe-yao}.

\begin{proposition} \label{prop-sigle-cycle}
	For $r\geq 2$, we have
	\begin{align*}
	\sum_{\Lc\in\oP([r])}\sum_{\tau\in\Cyc(\Lc)} (-1)^{|\Lc|-1}\int_{X^\Lc}  f(x_{\Lc(1)},x_{\Lc(2)},\dots,x_{\Lc(r)} ) \prod_{A\in\Lc}K_N (x_A,x_{\tau(A)}) \,\prod_{A\in\Lc}\omega(x_A) \\
	=\sum_{\eta:[r]\to[r]\setminus\{1\}} c_{\eta} \int_{X^r} \mathsf K_{N,r} \sum_{B\subset [r]}(-1)^{r-|B|} f\big( x^{B ,\eta}_1,x^{B ,\eta}_2,\dots,x^{B ,\eta}_r  \big) \prod_{\ell=1}^r \omega(x_\ell),
	\end{align*}
	where
	$$ \mathsf K_{N,r}:=\prod_{1\leq j\leq r} K_N(x_j,x_{j+1})   \quad\text{with}\quad x_{r+1}=x_1,   $$
	$$c_\eta:= (-1)^{ |\eta([r])| -1} \big|\eta([r])\big|^{-1} \binom{r}{\big|\eta([r])\big|} ^{-1},$$
	$$  x^{B ,\eta}_j:= x_{\eta(j)} \quad\text{for }\, j\in B, \qquad x^{B ,\eta}_j:=x_1 \quad\text{for }\, j\notin B. $$
\end{proposition}

\begin{proof}
	We start with a map $\eta:[r]\to[r]$ and let $\eta([r])=\{q_1,q_2,\dots,q_t\}$ with $q_1<q_2<\dots<q_t$. Put $A_j=\eta^{-1}(q_j)$. The reproducing property of $K_N$ in Lemma \ref{lem-reproducing} gives
	\begin{align}\label{eqknr-1}
	\int_{ X^{[r]\setminus   \eta([r])}} \mathsf K_{N,r} \,\prod_{\ell\in [r]\setminus   \eta([r])} \omega(x_\ell)= \prod_{\ell=1}^t K_N( x_{q_\ell},x_{q_{\ell+1}})   \quad\text{with}\quad x_{q_{t+1}}=x_{q_1}.
	\end{align}
	
	Note that $\Lc:=\{A_j:\,1\leq j\leq t\}$ is a partition of $[r]$, and we  define  $\tau\in \Sym(\Lc)$ by the single cycle $A_1\mapsto A_2\mapsto\cdots\mapsto A_t\mapsto A_1$. Adding $f$ and $\omega(x_\ell)$ for $\ell\in \eta([r])$ to the integrand in \eqref{eqknr-1}, we get
	\begin{align*}
	\int_{ X^{r } }\mathsf K_{N,r}  f \big( x_{\eta(1)},x_{\eta(2)},&\dots,x_{\eta(r)}  \big) \,\prod_{\ell=1}^r \omega(x_\ell) =\\
	&\int_{X^\Lc} f \big( x_{\Lc(1)},x_{\Lc(2)},\dots,x_{\Lc(r)}  \big)   \prod_{A\in\Lc}K_N(x_A,x_{\tau(A)}) \,\prod_{A\in\Lc} \omega(x_A).
	\end{align*}
	
	Now we count how many different $\eta$'s give the same $(\Lc,\tau)$. In the first step, we need to fix the image of $\eta$, having $\binom{r}{t}$ choices.
	After that, we assign a value to $\eta(1)$, having $t$ choices. At that moment, in fact, all the values of $\eta$ are already given by $\tau$ and $\eta(1)$ because $\tau$ is a single cycle. Thus, we have $t \binom{r}{t}=(-1)^{t-1}c_\eta^{-1}$ choices in total for every fixed  $(\Lc,\tau)$. Therefore,
	we conclude that 
	\begin{align}
	\sum_{\Lc\in\oP([r])}\sum_{\tau\in\Cyc(\Lc)} (-1)^{|\Lc|-1}\int_{X^\Lc}  f(x_{\Lc(1)},x_{\Lc(2)},\dots,x_{\Lc(r)} ) \prod_{A\in\Lc}K_N (x_A,x_{\tau(A)}) \,\prod_{A\in\Lc}\omega(x_A) \nonumber\\
	=\sum_{\eta:[r]\to[r]}c_\eta\int_{X^r}  \mathsf K_{N,r}  f\big( x_{\eta(1)},x_{\eta(2)},\dots,x_{\eta(r)}  \big) \,\prod_{j=1}^r \omega(x_j).\label{eqknr-2}
	\end{align} 
	
	\medskip
	\noindent
	\textbf{Claim:} We have
	\begin{align}\label{eqknr-3}
	\sum_{\eta:[r]\to[r]} c_\eta f\big( x_{\eta(1)},x_{\eta(2)},\dots,x_{\eta(r)}  \big) =\sum_{\eta:[r]\to[r]} c_\eta  \sum_{B\subset [r]}(-1)^{r-|B|} f\big( x^{B ,\eta}_1,x^{B ,\eta}_2,\dots,x^{B ,\eta}_r  \big).  
	\end{align}
	
	\begin{proof}[Proof of Claim]
		Notice that the left hand side of \eqref{eqknr-3} is the term in the right hand side when $B=[r]$. So we need to show that the terms for $B\neq [r]$ cancel with each other.  We take a proper subset $B$ of $[r]$. Fix $j_0\in [r]\setminus B$. For any map $\widetilde\eta:[r]\setminus\{j_0\}\to[r]$, we define the extensions $\eta_\ell$ for each $\ell=1,2,\dots,r$ by
		$$ \eta_\ell (j):=\widetilde \eta(j) \quad\text{for }\, j\neq j_0  \quad\text{and}\quad  \eta_\ell(j_0):=\ell. $$
		Observe that for all $\ell$,
		$$f\big( x^{B ,\eta_\ell}_1,x^{B ,\eta_\ell}_2,\dots,x^{B ,\eta_\ell}_r  \big) = f\big( x^{B ,\eta_1}_1,x^{B ,\eta_1}_2,\dots,x^{B ,\eta_1}_r  \big).  $$
		
		On the other hand, let $t:=\big|\widetilde\eta( [r]\setminus\{j_0\})\big|$. Then there are $t$ among the $\ell$'s such that $\big|\eta_\ell([r])\big|=t$ and $r-t$ number of $\ell$'s such that  $\big|\eta_\ell([r])\big|=t+1$. It follows that
		$$  \sum_{1\leq\ell\leq r} c_{\eta_\ell}= t (-1)^{t-1} t^{-1} \binom{r}{t}^{-1} +(r-t) (-1)^{t} (t+1)^{-1} \binom{r}{t+1}^{-1} =0.    $$
		This proves the Claim.
	\end{proof}
	
	For a map $\eta:[r]\to[r]$, if $s\in \eta^{-1}(1)$, then for any set $B\subset [r]$ not containing $s$, by the definition of $x^{B ,\eta}_s$, we have
	$$ f\big( x^{B,\eta}_1,x^{B,\eta}_2,\dots,x^{B,\eta}_r  \big)  = f\big( x^{B\cup\{s\},\eta}_1,x^{B\cup\{s\},\eta}_2,\dots,x^{B\cup\{s\},\eta}_r  \big).  $$
	Moreover, the coefficient $(-1)^{r-|B|}$ cancels with $(-1)^{r-|B\cup\{s\}|}$. This shows that  the terms indexed by these maps cancel in the sum of the right hand side of \eqref{eqknr-3}. Combining with \eqref{eqknr-2}, we finish the proof of the proposition.
\end{proof}

In the following, we will make a correspondence between the two different summation formulas in Proposition \ref{prop-cum} and Proposition \ref{prop-sigle-cycle}.

Fix a partition $\Fc:=\{ F_1,F_2,\dots,F_m  \}$ of $[2k]$. Define the two sets
\begin{align*}
\oA_\Fc:=\big\{ (\Qc,\sigma):\, \Qc \in\oP([2k]), \,\sigma\in\Sym(\Qc),\, \Qc\otimes \sigma=\Fc  \big\}
\end{align*}
and 
\begin{align*}
\oB_\Fc:= \Big\{ \big( (\Lc_1,\tau_1 ) ,(\Lc_2,\tau_2),\dots,(\Lc_m,\tau_m) \big):\,   \Lc_j\in \oP(F_j),\,\tau_j\in\Cyc(\Lc_j)           \Big\}.
\end{align*}

  Take a $(\Qc,\sigma)\in\oA_\Fc$. The condition $\Qc\otimes \sigma=\Fc$ implies that $\sigma$ is invariant on each $F_j$, i.e., for all $A\in \Qc$ with $A\subset F_j$, one has $\sigma(A)\subset F_j$. Thus, we can define $$\Lc_j:=\{ A\in\Qc:\, A\subset F_j\},\quad  \tau_j:=\sigma|_{\Lc_j}.$$ Clearly, $\tau_j$ is a single cycle and hence it gives an element in $\oB_\Fc$.
	Hence we have obtained a map from $\oA_\Fc\to\oB_\Fc$:
	\begin{equation}\label{mapoaob} (\Qc,\sigma) \mapsto \big( (\Lc_1,\tau_1),  (\Lc_2,\tau_2) ,\dots, (\Lc_m,\tau_m)     \big).   
		\end{equation}
It is not hard to check that   the map is
bijective, and the signs satisfy the following identity
\begin{align}\label{eq-sign}
\prod_{1\leq j\leq m} \sgn(\tau_j)=\sgn(\sigma).
\end{align}

\section{Proof of the main theorem}

Theorem \ref{thm-var} says that the standard deviation of $\Ec_N$ is
of order $\sqrt N$. To prove the Gaussian limit in
Theorem \ref{thm-main}, it therefore suffices to show that
$\kappa_k(\Ec_N)=o(N^{k/2})$ for every fixed $k\geq 3$.

The cumulant expansion in Proposition \ref{prop-cum} leaves only
connected configurations. Applying the single-cycle identity
of Proposition \ref{prop-sigle-cycle} rewrites  each Green function factor
as finite differences. The zero integral normalization of $G$
cancels the leading contribution from singleton blocks in the partition, and the
off-diagonal decay of the Bergman kernel controls the variables
within each remaining block of the partition.

\begin{theorem}\label{thm-cum}
	For $k\geq 3$ and  $N\geq 2k$, we have
	$$ \big|\kappa_k(\Ec_N)\big|\leq C_k N\log^k N $$
	for some constant $C_k>0$ independent of $N$.
\end{theorem}

The proof of Theorem \ref{thm-cum} is a bit long. We introduce some notation used in this section. For every partition $\Fc:=\{ F_1,F_2,\dots,F_m  \}$ of $[2k]$ satisfying $\Fc\vee \Rc_k=\{[2k]\}$, let $$S_\Fc :=\big\{ s\in [2k]:\, \{s\}\in \Fc\big\}$$  be the singletons  in $\Fc$, and list other elements in $\Fc$ as $B_1,B_2,\dots, B_t$ with at least two elements.  Let $R_j:=\{2j-1,2j\}$ and $\Rc_k$ be the partition $\{R_j:\, j\in [k]\}$ of $[2k]$.

\begin{lemma} \label{lem-singletone}
	There is no $j\in[k]$ such that $R_j\subset S_\Fc$ and $S_\Fc\neq [2k]$. In particular,  $t\geq 1$.
\end{lemma}

\begin{proof}
	If $R_j\subset  S_\Fc$ for some $j$, by definition of $\vee$, $R_j \in \Fc\vee \Rc_k$. We get a contradiction and prove the first assertion.
	
	If $S_\Fc= [2k]$, every element in the partition  $\Fc$ is singleton. Then $\Fc\vee \Rc_k=\Rc_k$. This is also a contradiction.
\end{proof}

Without loss of generality, we assume $F_1,F_2,\dots,F_t$ contain at least two elements and let $F_{t+1},\dots,F_m$ be the  singleton blocks.   
We list the elements in $F_\ell=\{u_{\ell,1},u_{\ell,2},\dots,u_{\ell,r_\ell}      \}$ for each $1\leq \ell\leq t$ and assume $u_{\ell,1}<u_{\ell,2}<\cdots<u_{\ell,r_\ell}$.  Put $u_{\ell,r_\ell+1}:=u_{\ell,1}$.  
Let $$J_\Fc:=\{j\in [k]:\, R_j\subset [2k]\setminus S_\Fc\}.$$  By Lemma \ref{lem-singletone}, $|J_\Fc|+|S_\Fc|=k$.
Define the function
	$$\varphi :[2k]\setminus S_\Fc \to [t] \quad\text{by}\quad        \varphi(s):=\ell \quad\text{if}\quad s\in F_\ell.  $$  
and the set of tuples of maps
$$ \mathcal V_\Fc:= \big\{ (\eta_1,\dots,\eta_t):\,  \eta_\ell \,\textit{ is a map }\,  F_\ell \to F_\ell\setminus\{u_{\ell,1}\}    \big\}. $$
For simplicity, we omit the dependence on $\Fc$ for some notations, as finally, we will only focus on one $\Fc$, which is enough for our proof.

We have the following formula for the cumulants of $\Ec_N$.

\begin{proposition} \label{prop-new-formula}
		For every $k\geq 2,N\geq 2k$, we have
	\begin{align*}
	\kappa_k(\Ec_N)=\sum_{\Fc\in \oP([2k])\atop \Fc\vee \Rc_k=\{[2k]\}}  \sum_{(\eta_1,\dots,\eta_t)\in\mathcal V_\Fc} \Big( \prod_{\ell=1}^t c_{\eta_\ell} \Big)\int_{X^{ 2k} } \prod_{\ell=1}^t \prod_{j=1}^{r_\ell} K_N( x_{u_{\ell, j}},x_{u_{\ell,j+1}} )  \prod_{q\in S_\Fc} K_N(x_q,x_q) \\   \prod_{j=1}^k \mathcal G^{[j]}_{\Fc,\eta_1,\dots,\eta_t}  \prod_{\ell=1}^t \prod_{j=1}^{r_\ell} \omega(x_{u_{\ell,j}}) \prod_{q\in S_\Fc} \omega (x_q), 
	\end{align*}
	where 
		\begin{align}
		\mathcal G^{[j]}_{\Fc,\eta_1,\dots,\eta_t}:=\tG\big(  &x_{\eta_{\varphi(2j-1)}(2j-1)},x_{\eta_{\varphi(2j)}(2j)}\big) -\tG\big(  x_{\eta_{\varphi(2j-1)}(2j-1)},x_{u_{\varphi(2j),1}}\big) \nonumber\\ 
	&	-\tG\big(  x_{u_{\varphi(2j-1),1}},x_{\eta_{\varphi(2j)}(2j)}\big) + \tG \big(x_{u_{\varphi(2j-1),1}},x_{u_{\varphi(2j),1}} \big)  \quad\text{if }\,  j \in J_\Fc, \label{nogj}
		\end{align} 
		\begin{align}\label{onegj}
		\mathcal G^{[j]}_{\Fc,\eta_1,\dots,\eta_t} : =\tG\big(  x_{2j-1},x_{\eta_{\varphi(2j)}(2j)}\big) -\tG\big(  x_{2j-1},x_{u_{\varphi(2j),1}}\big)\quad\text{if } \, 2j-1\in S_\Fc 
		\end{align}
		\begin{align*}
	\mathcal G^{[j]}_{\Fc,\eta_1,\dots,\eta_t}:  =\tG\big( x_{\eta_{\varphi(2j-1)}(2j-1)},x_{2j}\big) -\tG\big( x_{u_{\varphi(2j-1),1}}, x_{2j}\big)  \quad\text{if } \, 2j\in S_\Fc. 
		\end{align*}
	
\end{proposition}

\begin{proof}
	By Proposition \ref{prop-cum}, we   need to sum  the integrals
	$$ \oJ(\Qc,\sigma):= \sgn(\sigma)  \int_{X^{\Qc}} \prod_{j=1}^k \tG\big(x_{\Qc(2j-1)},x_{\Qc(2j)}\big)  \prod_{A\in\Qc} K_N(x_A ,x_{\sigma(A)})\,\prod_{A\in\Qc}\omega(x_A)   $$ 
	over all $(\Qc,\sigma)$ satisfying $(\Qc\otimes \sigma)\vee \Rc_k=\{[2k]\}$, where $\Qc\in\oP([2k])$ and $\sigma\in\Sym(\Qc)$. We can reorder  the sum, getting
	\begin{align*}
	\kappa_k(\Ec_N)=\sum_{\Qc\in \oP([2k])}\sum_{\sigma\in \Sym(\Qc)\atop(\Qc\otimes \sigma)\vee \Rc_k=\{[2k]\}} \oJ(\Qc,\sigma)=\sum_{\Fc\in\oP([2k])\atop \Fc\vee \Rc_k={\{[2k]\}}} \sum_{\Qc\in\oP([2k]) , \sigma\in\Sym(\Qc)\atop \Qc\otimes \sigma=\Fc}\oJ(\Qc,\sigma).
	\end{align*}
	Recall the bijective map $\oA_\Fc\to\oB_\Fc$ given in \eqref{mapoaob}.
	For each $\Fc$, we  use \eqref{eq-sign} to rewrite 
\begin{align*}
\oJ(\Qc,\sigma)=\prod_{j=1}^m \sgn(\tau_j)   \int_{X^{\Qc}} \prod_{j=1}^k \tG\big(x_{\Qc(2j-1)},x_{\Qc(2j)}\big)   \prod_{j=1}^m\prod_{A\in\Lc_j} K_N(x_A ,x_{\tau_j(A)})\,\prod_{j=1}^m\prod_{A\in\Lc_j}\omega(x_A).
\end{align*}
	Since \eqref{mapoaob} is bijective, the summation of $\oJ(\Qc,\sigma)$ over all $(\Qc,\sigma)\in\oA$ can be written as the summation over $\oB$ as
	\begin{align}
	\sum_{    \Lc_j\in \oP(F_j),\,\tau_j\in\Cyc(\Lc_j)}\prod_{j=1}^m  (-1)^{|\Lc_j|-1}   \int_{X^{\Lc}} \prod_{j=1}^k \tG\big(x_{\Lc(2j-1)},x_{\Lc(2j)}\big)   \prod_{j=1}^m\prod_{A\in\Lc_j} K_N(x_A ,x_{\tau_j(A)}) \nonumber\\
	\label{sum-lp}\,\prod_{j=1}^m\prod_{A\in\Lc_j}\omega(x_A),
	\end{align}
	where $\Lc:=\cup \Lc_j$.  We are going to apply Proposition \ref{prop-sigle-cycle}. Before that, we take out all the singletons $F_{t+1}\dots,F_m$.  In this case, $\tau_j=\id$ and $K_N(x_A,x_{\tau_j(A)})=K_N(x_A,x_A)$, which gives the factor $ \prod_{q\in S_\Fc} K_N(x_q,x_q)$.
	
	By Fubini's theorem, we can do the integration in the order of $x_s,s\in F_\ell$ for $1\leq \ell\leq t$. More precisely, we integrate the variables $x_{u_{\ell,1}},x_{u_{\ell,2}},\dots,x_{u_{\ell,r_\ell}}$ at the $\ell$-th step. 
Thus, applying Proposition \ref{prop-sigle-cycle} inductively for $t$ times, with
	$r_\ell, \Lc_\ell, \tau_\ell$ instead of $r,\Lc,\tau$, we obtain that \eqref{sum-lp} is 
\begin{align*}
  \sum_{(\eta_1,\dots,\eta_t)\in\mathcal V_\Fc} \Big( \prod_{\ell=1}^t c_{\eta_\ell} \Big)\int_{X^{2k} } \prod_{\ell=1}^t \prod_{j=1}^{r_\ell} K_N( x_{u_{\ell,j}},x_{u_{\ell,j+1}} )  \prod_{q\in S_\Fc} K_N(x_q,x_q) \\ \Phi_{\Fc,\eta_1,\dots,\eta_t}\,  \prod_{\ell=1}^t \prod_{j=1}^{r_\ell} \omega(x_{u_{\ell,j}}) \prod_{q\in S_\Fc} \omega (x_q)  
\end{align*}
for some function  $\Phi_{\Fc,\eta_1,\dots,\eta_t}$ on $X^{2k}$. Here, we let $x_{u_{\ell,r_\ell+1}}:=x_{u_{\ell,1}}$.
 
It remains to determine the function $\Phi_{\Fc,\eta_1,\dots,\eta_t}$.   Fix a vector $(\eta_1,\dots,\eta_t)\in\mathcal V_\Fc$. We apply   Proposition \ref{prop-sigle-cycle}  inductively on $\ell$ for $t$ times.  
   In $\ell$-th step, we  apply Proposition \ref{prop-sigle-cycle} without affecting $x_{u_{\ell',j}}$'s for $\ell'\neq \ell$. Moreover, the map 
$  (B_1,B_2,\dots,B_t)\mapsto \cup_{1\leq \ell\leq t} B_\ell       $
 gives a bijection from Cartesian product of subsets of $F_\ell$'s to the subsets of $[2k]\setminus S_\Fc$. Therefore,
 $$\Phi_{\Fc,\eta_1,\dots,\eta_t}= \sum_{B\subset [2k]\setminus S_\Fc}  (-1)^{ 2k-|S_\Fc| -|B|}  \prod_{j=1}^k \tG\big(y^{B}_{2j-1},y^{B}_{2j} \big)     $$
 where for each $1\leq s\leq 2k$,
 $$ y^B_s:=x_s \,\text{ if }\, s\in S_\Fc,\quad  y^B_s=x_{\eta_{\varphi(s)}(s)}  \,\text{ if }\, s\in B,\quad  y^B_s=x_{u_{\varphi(s),1}} \,\text{ if }\, s\notin S_\Fc \cup B.      $$

 For $s\in R_j$, the condition $s\in B$ is equivalent to $s\in A_j:=B\cap R_j$. It follows that $y^B_s=y^{A_j}_s$.
 Using that $  2k-|S_\Fc| -|B|=\sum_{1\leq j\leq k} |R_j\setminus S_\Fc|-|A_j|$, we have 
\begin{align*}
\Phi_{\Fc,\eta_1,\dots,\eta_t}&= \sum_{ B\subset [2k]\setminus S_\Fc , A_j\subset R_j\setminus S_\Fc \atop \cup A_j=B}  (-1)^{ 2k-|S_\Fc| -|B|}  \prod_{j=1}^k \tG\big(y^{A_j}_{2j-1},y^{A_j}_{2j} \big) \\
&=\prod_{j=1}^k \Big (   \sum_{A_j \subset R_j\setminus S_\Fc} (-1)^{|R_j\setminus S_\Fc|-|A_j| }   \tG\big(y^{A_j}_{2j-1},y^{A_j}_{2j} \big)\Big).
\end{align*}

We now evaluate the factors in the last product. By Lemma \ref{lem-singletone}, either both $2j-1,2j$ are not in $S_\Fc$ or exactly one of them are in $S_\Fc$. In the first case, $R_j\setminus S_\Fc=R_j$ and the  $A_j$'s we can take are $R_j, \{2j\},\{2j-1\}, \varnothing$. 
Then
\begin{align*}
\sum_{A_j \subset R_j \setminus S_\Fc}  (-1)^{|R_j|-|A_j| }   \tG \big(y^{A_j}_{2j-1},y^{A_j}_{2j} \big)   =\tG\big(  x_{\eta_{\varphi(2j-1)}(2j-1)},x_{\eta_{\varphi(2j)}(2j)}\big) -\tG\big(  x_{\eta_{\varphi(2j-1)}(2j-1)},x_{u_{\varphi(2j),1}}\big) \\
-\tG\big(  x_{u_{\varphi(2j-1),1}},x_{\eta_{\varphi(2j)}(2j)}\big) + \tG \big(x_{u_{\varphi(2j-1),1}},x_{u_{\varphi(2j),1}} \big). 
\end{align*}
In the second case, when $(2j-1)\in S_\Fc$. There are two $A_j$'s we can take: $\{2j\}$ and $\varnothing$. Then
\begin{align*}
\sum_{A_j \subset R_j\setminus S_\Fc}  (-1)^{1-|A_j| }   \tG \big(y^{A_j}_{2j-1},y^{A_j}_{2j} \big)   =\tG\big(  x_{2j-1},x_{\eta_{\varphi(2j)}(2j)}\big) -\tG\big(  x_{2j-1},x_{u_{\varphi(2j),1}}\big). 
\end{align*}
Similarly, when $2j \in S_\Fc$, 
 \begin{align*}
 \sum_{A_j \subset R_j\setminus S_\Fc}  (-1)^{1-|A_j| }   \tG \big(y^{A_j}_{2j-1},y^{A_j}_{2j} \big)   =\tG\big( x_{\eta_{\varphi(2j-1)}(2j-1)},x_{2j}\big) -\tG\big( x_{u_{\varphi(2j-1),1}}, x_{2j}\big). 
 \end{align*}
  We complete the proof of the proposition.
\end{proof}

In the summation of the formula in Proposition \ref{prop-new-formula}, for each fixed $\Fc$ and $(\eta_1,\dots,\eta_t)\in\mathcal V_\Fc$, the naive bound for the integral     is $O(N^{t}(\log N)^{|J_\Fc|})$ by using Lemma \ref{lem-kxx}, Proposition \ref{bergman-formula} and Lemma \ref{lem-ywlog}. In the following, we will select $t-1$ indices from $J_\Fc$ to create an $O_k(N^{1-t})$ factor.

\begin{lemma}\label{lem-reorder}
	One can re-order the sets $F_1,F_2,\dots,F_t$  and find distinct integers $\ee_\ell\in J_\Fc,1\leq \ell\leq t-1$  such that for all $1\leq \ell\leq t-1$,
	\begin{align*}
 	\varphi(2\ee_\ell-1)\neq \varphi(2\ee_\ell) ,\quad  \min\big( \varphi(2\ee_\ell-1),\varphi(2\ee_\ell)\big)=\ell.
	\end{align*}
\end{lemma}

\begin{proof}
	Take any set in $\Fc$ with at least two elements, denoted as $F_t$. Now we define $F_{t-1},\dots,F_1$ inductively.  Suppose $F_t,F_{t-1},\dots, F_{\ell+1}$ and $\ee_{t-1},\dots,\ee_{\ell+1} $ have been chosen, where $\ell\geq 1$. We are going to choose $F_\ell$ in $\Fc$. Let $B:=\cup_{\ell+1\leq q\leq t}F_q$. 
	
We first show	that
	 there exists some $j\in[k]$ such that $2j-1,2j\notin S_\Fc$ and $|R_j\cap B|=1$.  Suppose not. Then  $|R_j\cap B|=1$ will imply that the other element in $R_j$ is in $S_\Fc$.  Consider $A:=\cup_{ R_j\cap B\neq \varnothing } R_j$. Then $B\subset A \subset B\cup S_\Fc \neq [2k]$ and $A\setminus B\subset S_\Fc$.  This means $A$ is also the union of elements in $\Fc$, which contradicts with $\Fc\vee\Rc_k=\{[2k]\}$.  
	 
	 Therefore, we take this $j$ as $\ee_\ell$. One of $2\ee_\ell-1,2\ee_\ell$ is in $B$ and the other is in some $F\in \Fc$ not appearing in $F_t,F_{t-1},\dots, F_{\ell+1}$, which we denote by $F_\ell$. This gives the
	 three required conditions.
\end{proof}

We denote the index set $$L_\Fc:=\{\ee_1,\ee_2, \dots, \ee_{t-1} \}$$ chosen in   Lemma \ref{lem-reorder}.  The purpose of   re-ordering is the following: notice that  $\mathcal G^{[\ee_\ell]}_{\Fc,\eta_1,\dots,\eta_t} $ in \eqref{nogj}  depends on the variables $x_u$'s for $u\in F_{\varphi(2\ee_\ell-1)}\cup F_{\varphi(2\ee_\ell)}$.
Lemma \ref{lem-reorder} allows us to treat it as constant when   integrating with $x_u$'s for $u\in \cup_{q=1}^{\ell-1}F_{q}$.

\medskip

For each $1\leq\ell \leq t,2\leq j\leq r_\ell$, we define
$$d_{\ell,j}:= \dist(x_{u_{\ell,1}},x_{u_{\ell,j}} )  \quad\text{and}\quad d_\ell:=\max_{2\leq j\leq r_\ell}  d_{\ell,j}    .$$

\begin{lemma}\label{lem-C-ck}
	There exist   constants $C,c_k>0$ independent of $N$ such that for  $1\leq\ell\leq t$,
	\begin{align*}
 \prod_{1\leq j\leq  r_\ell} \big|K_N( x_{u_{\ell, j}},x_{u_{\ell,j+1}} ) \big| \leq C N^{r_\ell}  e^{-c_k \sqrt N d_\ell}  \prod_{2\leq j\leq  r_\ell} e^{-c_k \sqrt N d_{\ell,j}}.
	\end{align*}
\end{lemma}

\begin{proof}
	Let $$D_\ell:=\sum_{1\leq  j\leq r_\ell} \dist(x_{u_{\ell,j}},x_{u_{\ell,j+1}} ).$$
	The triangle inequality gives
	$$ d_\ell \leq D_\ell  \quad\text{and}\quad   \sum_{2\leq j\leq r_\ell} d_{\ell,j}\leq (r_\ell-1) D_\ell, $$
	which implies
	$ d_\ell+ \sum_{2\leq j\leq r_\ell} d_{\ell,j}\leq  r_\ell  D_\ell \leq 2k D_\ell $.
	Together with Lemma \ref{lem-globalber}, we get
	$$ \prod_{1\leq j\leq  r_\ell} \big|K_N( x_{u_{\ell, j}},x_{u_{\ell,j+1}} ) \big| \leq C N^{r_\ell}e^{-c \sqrt N D_\ell} \leq CN^{r_\ell} e^{ -c\sqrt N (2k)^{-1} ( d_\ell+\sum_{2\leq j\leq r_\ell} d_{\ell,j}) },  $$
	proving the desired inequality.
\end{proof}

For $j\in J_\Fc$, put
$$  v_{2j-1}:=   \eta_{\varphi(2j-1)}(2j-1), \quad v_{2j}:=\eta_{\varphi(2j)}(2j).  $$
 Fix a large  $a_k>0$ independent of $N$.
Let $ \Omega_{k,N} $  be the subset of $X^{2k-|S_\Fc|}$ where for all $j\in J_\Fc$,
\begin{align}\label{four-dist}
  \dist   \big(  x_{v_{2j-1}},x_{v_{2j}}\big) ,\dist\big(  x_{v_{2j-1}},x_{u_{\varphi(2j),1}}\big) 
,\dist\big(  x_{u_{\varphi(2j-1),1}},x_{v_{2j}}\big), \dist \big(x_{u_{\varphi(2j-1),1}},x_{u_{\varphi(2j),1}} \big)   
\end{align}
are all greater than $N^{-a_k}$, whenever they are not identically zero (this happens when the two variables coincide). Immediately, from the definition of  $\mathcal G^{[j]}_{\Fc,\eta_1,\dots,\eta_t}$ for $j\in J_\Fc$ in \eqref{nogj} and using \eqref{decom-green}, we see that   (recall $\tG(x,x)=0$)
\begin{align}\label{Gbound1}
\big|\mathcal G^{[j]}_{\Fc,\eta_1,\dots,\eta_t}\big| \lesssim -\log N^{-a_k}=a_k\log N  \quad\text{on}\quad \Omega_{k,N}.
\end{align}

For the integral of $\prod \big| \mathcal G^{[j]}_{\Fc,\eta_1,\dots,\eta_t} \big|$ outside $\Omega_{k,N}$ we have the following easy bound.

\begin{lemma}\label{lem-omegaout}
	We have
	$$\int_{X^{ 2k-|S_\Fc|} \setminus \Omega_{k,N}} \prod_{\ell=1}^t \prod_{j=1}^{r_\ell} \big| K_N( x_{u_{\ell, j}},x_{u_{\ell,j+1}} )  \big|   \prod_{j\in J_\Fc} \big| \mathcal G^{[j]}_{\Fc,\eta_1,\dots,\eta_t} \big|\, \prod_{\ell=1}^t \prod_{j=1}^{r_\ell} \omega(x_{u_{\ell,j}})\lesssim_k N^{-2}.$$
\end{lemma}

\begin{proof}
	Take a $j\in J_\Fc$. We have $v_{2j-1}\neq u_{\varphi(2j),1}$.
	Since $|K_N(x,y)|\lesssim N$ and  $\prod |\mathcal G^{[j]}_{\Fc,\eta_1,\dots,\eta_t}|^2$ is integrable, by Cauchy-Schwarz inequality,
	the left hand side is
	\begin{align*}
	&\lesssim_k N^{2k-|S_\Fc|} \omega^{2k-|S_\Fc|}\big(X^{2k-|S_\Fc|}\setminus\Omega_{k,N} \big)^{1/2} \\
	&\leq N^{2k-|S_\Fc|} \omega^{2k-|S_\Fc|}\Big(\{ \dist (  x_{v_{2j-1}},x_{u_{\varphi(2j),1}} )\leq N^{-a_k} \}  \times X^{2k-|S_\Fc|-1}  \Big)^{1/2}
	\\&= N^{2k-|S_\Fc|} \omega\Big(\{ \dist (  x_{v_{2j-1}},x_{v_{2j}} )\leq N^{-a_k} \}   \Big)^{1/2}\lesssim N^{2k-|S_\Fc|} N^{-a_k}.
	\end{align*} 
	 By taking  $a_k\geq 2k+2$, we get the desired inequality.
\end{proof}

To bound the integral inside $\Omega_{k,N}$ is not trivial.
Let $$ \mathsf D_j:= \dist (x_{u_{\varphi(2j-1),1}},x_{u_{\varphi(2j),1}}  ).$$
We have the following more delicate bound for $\mathcal G^{[j]}_{\Fc,\eta_1,\dots,\eta_t}$.

\begin{lemma}\label{lem-Gbound}
	There exists  a constant $C >0$ independent of $k,N$ such that 
  $$  \big |\mathcal G^{[j]}_{\Fc,\eta_1,\dots,\eta_t} \big|\leq C \big(d_{\varphi(2j-1)}^2+ d_{\varphi(2j)}^2 \big) \mathsf D_j^{-2}  \quad\text{on} \quad\Omega_{k,N} $$
  when $\mathsf D_j >9 (d_{\varphi(2j-1)}+d_{\varphi(2j)})$. 
\end{lemma}

\begin{proof}
	The assumption on $\mathsf D_j$ gives $\varphi(2j-1)\neq \varphi(2j)$ and hence  all the quantities in \eqref{four-dist} are non-zero.
	By definition, 
	$$\dist\big(  x_{u_{\varphi(2j-1),1}},x_{v_{2j-1}}\big)\leq  d_{\varphi(2j-1)} \quad\text{and}    \quad   \dist\big(  x_{u_{\varphi(2j),1}},x_{v_{2j}}\big)\leq  d_{\varphi(2j)}.$$
	It follows from $\mathsf D_j > 9 (d_{\varphi(2j-1)}+d_{\varphi(2j)})$ that all the quantities in \eqref{four-dist} are $\geq \mathsf D_j/2$. Thus, $\big |\mathcal G^{[j]}_{\Fc,\eta_1,\dots,\eta_t} \big|\lesssim |\log\mathsf D_j|$ and	it is enough to assume   $\mathsf D_j$ is small. In this case, all the variables defining $\mathcal G^{[j]}_{\Fc,\eta_1,\dots,\eta_t}$ are in the local normal coordinate centering at $x_{u_{\varphi(2j-1),1}}$.  Using  \eqref{decom-green} and that $\varrho$ is smooth on this coordinate, it is enough to show
	\begin{align}\label{xi-zeta}
	 \big|   \log|\xi-\zeta |-\log|z-\zeta|-\log|\xi|     +\log|z|  \big|\lesssim \big(d_{\varphi(2j-1)}^2+ d_{\varphi(2j)}^2 \big) \mathsf D_j^{-2}  
	\end{align}
	where $z,\xi,\zeta$ are the coordinates of $x_{u_{\varphi(2j),1}},x_{v_{2j}},x_{v_{2j-1}}$ respectively. Since the distance in the normal coordinate is comparable with the distance in $X$, we have
	\begin{align*} 
	|\zeta|\leq 1.1 d_{\varphi(2j-1)}, \quad |\xi-z|\leq  1.1 d_{\varphi(2j)}, \quad |z|\geq 0.9 \mathsf D_j. 
	\end{align*}  
	Then it is not hard to see that for $w\in [\xi,z]$,
	$$ |w|\geq 0.77 \mathsf D_j \quad\text{and}\quad |w-\zeta|\geq 0.64 \mathsf D_j,  $$

	Consider the function $f(w):=\log|w-\zeta|-\log|w|$. By mean value theorem,
	$$|f(\xi)-f(z)|\leq |\xi-z| \sup_{w\in[\xi,z]} |\nabla f(w)|=|\xi-z| \sup_{w\in[\xi,z]} \Big|{w-\zeta \over |w-\zeta|^2 } -{w\over |w|^2} \Big|\lesssim |\xi-z|{|\zeta|\over \mathsf D_j^2}.   $$
	  This gives \eqref{xi-zeta} and the lemma follows.
\end{proof}

\begin{lemma}\label{lem-gjglobal}
	For $j\in L_\Fc$ and
  $c_k>0$ independent of $N$, we have on $\Omega_{k,N}$,
	\begin{align*}
	  \big |\mathcal G^{[j]}_{\Fc,\eta_1,\dots,\eta_t} \big| e^{-c_k \sqrt N (d_{\varphi(2j-1)}+d_{\varphi(2j)})} \lesssim_k \big(\log N  e^{-c_k \sqrt N \mathsf D_j/9}+N^{-1}\mathsf D_j^{-2}\big).
	\end{align*}
\end{lemma}

\begin{proof}
	For $j\in L_\Fc$, Lemma \ref{lem-reorder} says that $\varphi(2j-1)\neq \varphi(2j)$ and hence $\mathsf D_j$ is not identically zero.
	When $\mathsf D_j \leq  9( d_{\varphi(2j-1)}+d_{\varphi(2j)})$, we use \eqref{Gbound1} to obtain the left hand side is
	$$ \lesssim_k  \log N  e^{-c_k \sqrt N \mathsf D_j/9}. $$
	When $\mathsf D_j > 9( d_{\varphi(2j-1)}+d_{\varphi(2j)})$,   Lemma \ref{lem-Gbound} gives the left hand side is 
	\begin{align*}
	\lesssim \big(d_{\varphi(2j-1)}^2+ d_{\varphi(2j)}^2 \big) \mathsf D_j^{-2}e^{-c_k \sqrt N (d_{\varphi(2j-1)}+d_{\varphi(2j)})} \lesssim_k  N^{-1}\mathsf D_j^{-2}.
	\end{align*}
	The lemma follows.
\end{proof}

Recall from Lemma \ref{lem-reorder} that $\mathsf D_j,j\in L_\Fc$ depend on $x_{u_{\ell,1}},1\leq \ell\leq t$. Define $\Gamma_t:=X$ and for $1\leq s\leq t-1$, define   
$$  \Gamma_s:=\big\{( x_{u_{s,1}}, x_{u_{s+1,1}},\dots,  x_{u_{t,1} }  )\in  X^{t-s+1}  :\,  \mathsf D_j>N^{-a_k} \,\text{ for }\, j=\ee_s,\ee_{s+1},\dots,\ee_{t-1}    \big\}.         $$

\begin{lemma}\label{lem-inductive}
	For $c_k'>0$ independent of $N$, we have
	\begin{align*}
\int_{\Gamma_1} \prod_{j\in L_\Fc} \big(\log N  e^{-c'_k \sqrt N \mathsf D_j }+N^{-1}\mathsf D_j^{-2}\big) \, \prod_{1\leq\ell\leq t} \omega(x_{u_{\ell,1}})\lesssim_k  N^{-t+1}\log^{t-1}N.
	\end{align*} 
\end{lemma}

\begin{proof}
	Lemma \ref{lem-reorder} says that  for each $1\leq s\leq t-1$, $\mathsf D_{\ee_s}=\dist(x_{u_{s,1}}, x_{u_{\ell',1}})$ for some $\ell'>s$ and $\mathsf D_{\ee_{\ell}}$ is independent of $x_{u_{s,1}}$ for all $\ell> s$.  Thus, by Fubini's theorem, and using that 
	\begin{align}\label{supx}
	\sup_{x\in X}\int  e^{-c_k' \sqrt N \dist(x,y)}\,\omega(y) \lesssim _k N^{-1},\quad  \sup_{x\in X}\int_{\dist(x,y)>N^{-a_k}} \dist(x,y)^{-2}\,\omega(y)\lesssim_k \log N,     
	\end{align}     
	we have
	\begin{align*}
	&\int_{\Gamma_{s}} \prod_{ \ell=s}^{ t-1}  \big(\log N  e^{-c'_k \sqrt N \mathsf D_{\ee_\ell} }+N^{-1}\mathsf D_{\ee_\ell}^{-2} \big) \, \prod_{\ell=s}^{ t} \omega(x_{u_{\ell,1}})\\
	&	=\int_{\Gamma_{s+1}} \prod_{ \ell=s+1}^{ t-1} \big (\log N  e^{-c'_k \sqrt N \mathsf D_{\ee_\ell} }+N^{-1}\mathsf D_{\ee_\ell}^{-2} \big) \\
	&\qquad \qquad  \int_{\mathsf D_{\ee_s}>N^{-a_k}} \big(\log N  e^{-c'_k \sqrt N \mathsf D_{\ee_s} }+N^{-1}\mathsf D_{\ee_s}^{-2} \big) \, \omega(x_{u_{s,1}}) \,   \prod_{ \ell=s+1}^{ t} \omega(x_{u_{\ell,1}})\\
	&\lesssim_k N^{-1}\log N \int_{\Gamma_{s+1}} \prod_{ \ell=s+1}^{ t-1} \big (\log N  e^{-c'_k \sqrt N \mathsf D_{\ee_\ell} }+N^{-1}\mathsf D_{\ee_\ell}^{-2} \big)\,   \prod_{ \ell=s+1}^{ t} \omega(x_{u_{\ell,1}}).
	\end{align*} 
	By induction, the integral in the lemma is $O_k( N^{-t+1}\log^{t-1}N )$. 
\end{proof}

We are ready to bound the integral of $\prod\big| \mathcal G^{[j]}_{\Fc,\eta_1,\dots,\eta_t} \big|$ on $\Omega_{k,N}$.

\begin{lemma}\label{lem-omegabound}
	There exists  a constant $C_k>0$ independent of $N$ such that
  \begin{align*}
  \int_{ \Omega_{k,N}} \prod_{\ell=1}^t \prod_{j=1}^{r_\ell} \big| K_N( x_{u_{\ell, j}},x_{u_{\ell,j+1}} )  \big|   \prod_{j\in J_\Fc} \big| \mathcal G^{[j]}_{\Fc,\eta_1,\dots,\eta_t} \big|\, \prod_{\ell=1}^t \prod_{j=1}^{r_\ell} \omega(x_{u_{\ell,j}})\leq C_k N (\log N)^{|J_\Fc|}.
  \end{align*}
\end{lemma}

\begin{proof}
	 Applying Lemma \ref{lem-C-ck} to the product of $K_N$'s and  using \eqref{Gbound1} to the $ \mathcal G^{[j]}_{\Fc,\eta_1,\dots,\eta_t}$'s for $j\in J_\Fc\setminus L_\Fc$,   the integrand is 
	\begin{align*}
	&\lesssim_k \prod_{1\leq\ell\leq t} \Big(C N^{r_\ell}  e^{-c_k \sqrt N d_\ell}  \prod_{2\leq j\leq  r_\ell} e^{-c_k \sqrt N d_{\ell,j}}\Big) \prod_{j\in J_\Fc\setminus L_\Fc} (\log N) \prod_{j\in   L_\Fc}\big| \mathcal G^{[j]}_{\Fc,\eta_1,\dots,\eta_t} \big|\\
	&\lesssim_k N^{2k-|S_\Fc|} \prod_{1\leq\ell\leq t}e^{-c_k \sqrt N  d_\ell}  \prod_{1\leq\ell\leq t \atop 2\leq j\leq  r_\ell }e^{-c_k \sqrt N d_{\ell,j}} (\log N)^{|J_\Fc|-|L_\Fc|}\prod_{j\in   L_\Fc}\big| \mathcal G^{[j]}_{\Fc,\eta_1,\dots,\eta_t} \big|.
	\end{align*}
	
	Since $|L_\Fc|=t-1$, and $\varphi(2j-1)\neq \varphi(2j)$ for $j\in L_\Fc$, we have
	\begin{align*}
	     \sum_{j\in L_\Fc} (d_{\varphi(2j-1)}+d_{\varphi(2j)})      \leq (t-1) \sum_{1\leq \ell \leq t} d_\ell \leq k\sum_{1\leq \ell \leq t}d_\ell.
	\end{align*}
	Hence by Lemma \ref{lem-gjglobal},
	\begin{align*}
	&\prod_{1\leq\ell\leq t}e^{-c_k \sqrt N  d_\ell}\prod_{j\in   L_\Fc}\big| \mathcal G^{[j]}_{\Fc,\eta_1,\dots,\eta_t} \big| \leq e^{-c_k \sqrt Nk^{-1}\sum_{j\in L_\Fc} (d_{\varphi(2j-1)}+d_{\varphi(2j)})}\prod_{j\in   L_\Fc}\big| \mathcal G^{[j]}_{\Fc,\eta_1,\dots,\eta_t} \big|\\
	&=\prod_{j\in   L_\Fc}    \Big(  \big| \mathcal G^{[j]}_{\Fc,\eta_1,\dots,\eta_t} \big| e^{-c_k \sqrt Nk^{-1} (d_{\varphi(2j-1)}+d_{\varphi(2j)})} \Big)
	\lesssim_k \prod_{j\in   L_\Fc} \big (\log N  e^{-c_k' \sqrt N \mathsf D_j }+N^{-1}\mathsf D_j^{-2}\big).
	\end{align*}
	We deduce that the integral in the lemma is 
	\begin{align}
	\lesssim_k N^{2k-|S_\Fc|} (\log N)^{|J_\Fc|-t+1}\int_{ \Omega_{k,N}} &\prod_{1\leq\ell\leq t \atop 2\leq j\leq  r_\ell }e^{-c_k \sqrt N d_{\ell,j}}  \nonumber\\
	&\prod_{j\in   L_\Fc}  \big(\log N  e^{-c'_k \sqrt N \mathsf D_j }+N^{-1}\mathsf D_j^{-2} \big) \, \prod_{1\leq\ell\leq t \atop 1\leq j\leq  r_\ell } \omega(x_{u_{\ell,j}}).\label{int-omega-1}
	\end{align}
	
	Since $\mathsf D_j$ only depends on $x_{u_{\varphi(2j-1),1}},x_{u_{\varphi(2j),1}}$, and moreover, $x_{u_{\ell,j}},2\leq j\leq r_\ell$ only appear once in $d_{\ell,j}$. By Fubini's theorem, we can integrate $x_{u_{\ell,j}}$ for $2\leq j\leq r_\ell$ separately.  Using \eqref{supx}, we see that the product of $e^{-c_k \sqrt N d_{\ell,j}}$'s will contribute $O_k(N^{-\sum_{1\leq \ell\leq t} (r_\ell-1)})=O_k( N^{-2k+|S_\Fc|+t})$ to the integral.  On the other hand, the projection of $\Omega_{k,N}$ to the variables $x_{u_{\ell,1}},1\leq \ell\leq t$ is contained in $\Gamma_1$.  Therefore, \eqref{int-omega-1} is 
	\begin{align*} 
	\lesssim_k N^{2k-|S_\Fc|} (\log N)^{|J_\Fc|-t+1} N^{-2k+|S_\Fc|+t} \int_{\Gamma_1} \prod_{j\in   L_\Fc}  \big(\log N  e^{-c'_k \sqrt N \mathsf D_j }+N^{-1}\mathsf D_j^{-2} \big) \, \prod_{1\leq\ell\leq t} \omega(x_{u_{\ell,1}}).
	\end{align*}
Together with Lemma \ref{lem-inductive}, 
we conclude   that  \eqref{int-omega-1} is
 $$ \lesssim_k  N^{2k-|S_\Fc|} (\log N)^{|J_\Fc|-t+1} N^{-2k+|S_\Fc|+t} N^{-t+1}\log^{t-1}N =N (\log N )^{|J_\Fc|}.$$
	The proof of the lemma is complete.
\end{proof}

Now we can complete the proof of the main result of this section.
 
\begin{proof}[End of the proof of Theorem \ref{thm-cum}]
	By Proposition  \ref{prop-new-formula}, we only need to bound   \begin{align}\label{pt-1}
	\bigg|\int_{X^{ 2k} } \prod_{\ell=1}^t \prod_{j=1}^{r_\ell} K_N( x_{u_{\ell, j}},x_{u_{\ell,j+1}} )  \prod_{q\in S_\Fc} K_N(x_q,x_q)    \prod_{j=1}^k \mathcal G^{[j]}_{\Fc,\eta_1,\dots,\eta_t}  \prod_{\ell=1}^t \prod_{j=1}^{r_\ell} \omega(x_{u_{\ell,j}}) \prod_{q\in S_\Fc} \omega (x_q)\bigg|
	\end{align}
	for   fixed $\Fc$ and $(\eta_1,\dots,\eta_t)$ because the total number of choices depends only  on $k$.

	First, we integrate the $x_q$'s   over the measure $\prod_{q\in S_\Fc} K_N(x_q,x_q)\,  \omega(x_q)$ for $q\in S_\Fc$. If $2j-1\in S_\Fc$, by Lemma \ref{lem-kxx}, using formula \eqref{onegj} and that $\int_X G(x,\cdot)\,\omega=0$, we have
	$$ \int_X \mathcal G^{[j]}_{\Fc,\eta_1,\dots,\eta_t}   K_N(x_{2j-1},x_{2j-1})\,  \omega(x_{2j-1}) = \int_X \mathcal G^{[j]}_{\Fc,\eta_1,\dots,\eta_t} O(1)\,  \omega(x_{2j-1})=O(1).$$
	 The case when $2j\in S_\Fc$ is analogous. Therefore, \eqref{pt-1} is
	\begin{align} \label{pt-2}
	\lesssim_k \int_{X^{ 2k-|S_\Fc|} } \prod_{\ell=1}^t \prod_{j=1}^{r_\ell} \big| K_N( x_{u_{\ell, j}},x_{u_{\ell,j+1}} )  \big|   \prod_{j\in J_\Fc} \big| \mathcal G^{[j]}_{\Fc,\eta_1,\dots,\eta_t} \big|\, \prod_{\ell=1}^t \prod_{j=1}^{r_\ell} \omega(x_{u_{\ell,j}}),
	\end{align}
    which by Lemmas \ref{lem-omegabound} and \ref{lem-omegaout}, is 
    $\lesssim_k  N (\log N)^{|J_\Fc|}$. Using the fact that $|J_\Fc|\leq k$, we finish the proof of Theorem \ref{thm-cum}.
\end{proof}

 We now write down the proof of our main theorem.
 
 \begin{proof}[Proof of Theorem \ref{thm-main}]
 	Let $$Y_N:=  { \Ec_N- \bE(\Ec_N) \over \sqrt{ \Var (\Ec_N)    } }. $$
 	Then $\bE(Y_N)=0$ and $\Var(Y_N)=1$ for all $N$. By Theorems \ref{thm-var} and \ref{thm-cum}, for $k\geq 3$ and $N$ large, we have
 	$$\big| \kappa_k (Y_N)\big|= \Var(\Ec_N)^{-k/2} \big|\kappa_k (\Ec_N)\big| \lesssim_k N^{-k/2} N\log^k N=N^{1-k/2}\log^k N.  $$
 	The right-hand side tends to zero as $N\to\infty$. This shows $Y_N\to \oN(0,1)$ in law.
 \end{proof}

\section{Remark on CLT for logarithmic energy on complex manifolds} \label{sec-remk}

In this section, let $X$ be a compact K\"ahler manifold  of dimension $d\geq 1$, endowed with a positive holomorphic line bundle $\oL$. Fix a Hermitian metric $\fh$ on $\oL$ such that the curvature form $c_1(\oL,\fh)$ is strictly positive, which induces a K\"ahler form  $\omega$ with $\int_X \omega^d=1$. By   Riemann-Roch-Hirzebruch theorems (c.f.\ \cite[Theorem 1.4.6]{ma-book})
$$  N:=\dim_{\mathbb{C}} H^0(X, \oL^n) = V_\oL n^d +O(n^{d-1}) \quad\text{where}\quad  V_\oL:=\frac{1}{d!}\int_X c_1(\oL,\fh)^d. $$

The Bergman kernel of $H^0 (X,\oL^n)$, still denoted by  $K_{N}(x,y)$, has the following   Tian–Catlin–Zelditch expansion \cite{cat-ber,tian-jdg,zel-imrn-sze}
$$  K_N(x,x)=V_\oL \big( n^d +\fb_1(x) n^{d-1}+\fb_2(x) n^{d-2}+\cdots +\fb_d(x) +O(n^{-1})  \big),  $$
where  $\fb_j$'s are smooth real-valued functions independent of $n$.

Consider  the determinantal point process given by
$$\rho^N(x_1,x_2,\dots,x_{N})= {1\over N !}   \det \big( K_N(x_j,x_\ell) \big)_{1\leq j,\ell\leq N} \,\omega^d(x_1)\omega^d(x_2)\cdots\omega^d(x_{N})$$ 
and the logarithmic energy 
$$  \Ec^{\log}_{N}(\mathbf x):= \sum_{1\leq j,k\leq N, j\neq k}\log\dist (x_j,x_k).   $$

The variance and cumulants of $\Ec^{\log}_{N}$ will depend on the smallest $j$ that makes the function 
$$H_j(x):= \int_X \log\dist(x,y) \fb_j (y)\,\omega^d(y)  \qquad\text{(set $\fb_0=\mathbf 1$)}   $$ 
to be non-constant, whenever such an index exists. 
Therefore, one needs an even longer computation than the Riemann surfaces case.

 The centered logarithmic energy, divided by its standard deviation, also converges in law to a standard Gaussian as $N\to\infty$, while the proof and the variance estimates will appear in our forthcoming work \cite{wu-clt-2}.

\medskip

\end{document}